\documentclass[11pt]{amsproc}

\usepackage[dvips]{graphicx} 
\usepackage[counterclockwise]{rotating}

\usepackage{mathrsfs}
\newcommand{\J}{\mathscr J}

\usepackage{hyperref}

\renewcommand{\baselinestretch}{1.1}

\title[Elliptic parameters]{%
Elliptic parameters and defining equations \\
for elliptic fibrations on a Kummer surface}
\date{June, 2006}

\author[M.~Kuwata]{Masato Kuwata}
\address[M.~Kuwata]{Faculty of Economics
Chuo University\\
742-1 Higashi-Nakano, Hachioji-shi\\
Tokyo 192-0393,  Japan}
\email{kuwata@tamacc.chuo-u.ac.jp}

\author[T.~Shioda]{Tetsuji Shioda}
\address[T.~Shioda]{Department of Mathematics
Rikkyo University\\
3-34-1 Nishi-Ikebukuro, Toshima-ku\\
Tokyo 171-8501, Japan}
\email{shioda@rkmath.rikkyo.ac.jp}


\newtheorem{thm}{Theorem} 

\newtheorem{lem}[thm]{Lemma}

\newtheorem{problem}[thm]{Problem} 
\theoremstyle{definition}
\newtheorem{example}{Example}[section]
\newtheorem*{remark}{Remark}
\newcommand{\A}{{\bf A}}

\newcommand{\Q}{{\bf Q}}
\newcommand{\Z}{{\bf Z}}

\newcommand{\C}{{\bf C}}
\newcommand{\Proj}{{\bf P}}
\newcommand{\Hom}{\operatorname{Hom}}
\newcommand{\Km}{\operatorname{Km}}

\renewcommand{\l}{\lambda}
\newcommand{\ds}{\displaystyle}
\makeatletter
 
  \@addtoreset{equation}{section}
\makeatother

\begin{document}

\begin{abstract} 
We pose the problem to determine explicit defining equations of various elliptic fibrations on a given $K3$ surface, and study the case of the Kummer surfaces of the product of two elliptic curves.
\end{abstract}

\maketitle

\section{Introduction} 
\subsection{Problem setting}

Let $X$ be a $K3$ surface defined over a base field $k$, and let $k(X)$
denote its function field. Suppose $f:X \to \Proj^1$ is an elliptic
fibration on $X$ with a section $O$. Then it defines a non-constant
function $u=f(x)$ $(x \in X)$, and hence an element $u \in k(X)$. We
call $u$ the {\em elliptic parameter} for the elliptic fibration $f$.
(Actually $u$ is unique only up to the linear fractional
transformations, but to fix the idea, we always choose one $u$. Note
that the subfield $k(u)$ of $k(X)$ is uniquely defined by $f$).

Now let $E$ denote the generic fiber of $f$. Then $E$ is an elliptic
curve defined over $k(u)$ such that the function field $k(u)(E)$ is
isomorphic to $k(X)$ as the extensions of $k$.

\begin{problem} Given a $K3$ surface $X/k$ and an elliptic fibration $f$,
determine (i)  the elliptic parameter $u$ for $f$, (ii)
 the defining equation of the elliptic curve $E /k(u)$, and (iii) the Mordell-Weil lattice (MWL) $E(k(u))$.
\end{problem}

\begin{problem} Given a $K3$ surface $X/k$,  determine all the (essentially distinct) elliptic parameters.
\end{problem}

 Problem 2 is  a combination of Problem 1 and
the following standard problem:
\begin{problem}
Given a $K3$ surface $X/k$,  classify the elliptic fibrations 
$f:X \to \Proj^1$ up to isomorphisms.
\end{problem}

\subsection{Main results}
In this paper, we focus on the case of Kummer surfaces $X=\Km (A)$,
where $A=C_1 \times C_2$ is a product of two elliptic curves, and
assume $k$ is an algebraically closed field of characteristic different
from~2.

In this case, Problem~3 has been solved by Oguiso \cite{O} under the
assumption
\begin{itemize}
\item[$(\#)$] $C_1, C_2$ are not isogenous to each other and $k = \C$
(the field of complex numbers).
\end{itemize}
Namely he classifies the configuration
of singular fibers on such a Kummer surface $X$ into eleven types
$\J_1, \ldots, \J_{11}$, and determines the number of the
isomorphism classes for each type.

Our main results can be stated as follows: we solve Problem~1 for each
type of Oguiso's list (without assuming $(\#)$), and thus solve
Problem~2 under the assumption $(\#)$. More details will be given in
\S1.5 and~1.7 after we fix the notation and review some known cases.

\subsection{Notation}

By a $(-2)$-curve we mean a smooth rational curve on $X$ whose 
self-intersection number is $-2$. (It is called a ``nodal curve'' in
Oguiso \cite{O}.) It is known (cf. \cite{Kodaira}) that all irreducible
components of a reducible fiber in an elliptic fibration are
$(-2)$-curves.

\begin{figure}[!ht]
\begin{center}
\includegraphics[scale=0.72]{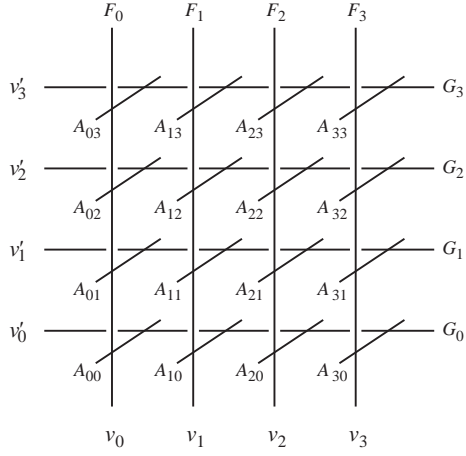}
\caption{double Kummer pencil}
\label{double_Kummer_pencil}
\end{center}
\end{figure}

We have a configuration of twenty-four $(-2)$-curves on X, called
the {\em double Kummer pencil} (see~Fig.~\ref{double_Kummer_pencil}, 
cf.~\cite{IS}). It consists of the 16
exceptional curves $A_{ij}$ arising from the minimal resolution $X \to
A/\iota_A$, plus the 8 curves $F_i, G_j$ obtained as the image of $v_i
\times C_2$ or $C_1\times v_j'$ under the rational map $A \to S$. Here
$\{v_i\}$ (or $\{v_i'\}$) denote the 2-torsion points of $C_1$
 (resp.~$C_2$) ($i,j \in I=\{0,1,2,3\}$), and $\iota_A$ denotes the inversion
automorphism of $A$. These curves will be referred to as the {\em basic
curves} below.

Suppose that the elliptic curve $C_i$ is defined by the Legendre form
\[
C_{i} : y_i^2=x_i (x_i -1)(x_i-\lambda_i) \qquad \l_{i}\neq 0,1.
\]
We order the 2-torsion
points by $v_1=(0,0), v_2=(1,0), v_3=(\lambda_1,0)$, with $v_0$
denoting the origin of $C_1$; similarly for $v_j'$ and $C_2$.

The function field $k(X)$ is equal to the subfield of the function field
$k(A)=k(x_1,y_1,x_2,y_2)$ consisting of the elements invariant under the
inversion $(x_1,y_1,x_2,y_2)\mapsto (x_1,-y_1,x_2,-y_2)$, namely we
have $$k(X)= k(x_1,x_2, t), \quad t= \frac{y_2}{y_1},$$ where $x_1,x_2$
and $t$ are naturally regarded as functions on $X$, satisfying the
relation
\begin{equation}\label{Kummer} 
x_{1} ( x_{1}-1 ) ( x_{1}-\l_{1}) t^{2} 
= x_{2} ( x_{2}-1 ) ( x_{2}-\l_{2} ).
\end{equation}

\subsection{Examples}
We start from the most classical and elementary example:

\begin{figure}[!ht]
\begin{center}
\[
\includegraphics[scale=0.72]{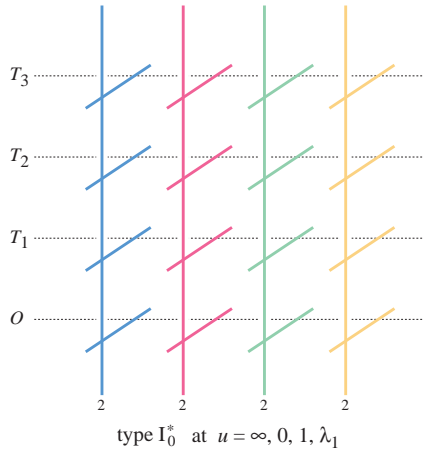}
\]
\caption{Kummer pencil (type $\J_{4}$)}
\label{J4}
\end{center}
\end{figure}
\begin{example}[Kummer pencils]  
The projection of $A$ to the first factor induces an elliptic
fibration $\pi_1: X \to \Proj^1$ with four singular fibers of type
I${}_0^{*}$:
$$  \Phi_i = 2 F_i + \sum_j{A_{ij}}  $$
(see Fig.~\ref{J4}).
This $\pi_1$ and the similar $\pi_2$ (obtained from the second
projection) are respectively called the first or second Kummer pencil
on $X$. The elliptic parameter for $\pi_1$ (or $\pi_2$) is obviously
given by the function $x_1$ (resp. $x_2$) in $k(X)$. (This belongs to
type $\J_4$ in \cite{O}, and $\pi_1$ and $\pi_2$ are the two
representatives of isomorphism classes, if $C_1,C_2$ are not
isogenous.)

The defining equation of the generic fiber over $k(x_1)$ is easily
obtained (see \S2.3), which is isomorphic to the constant curve $C_2$
over the quadratic extension $k(x_1,y_1)=k(C_1)$ of $k(x_1)$. The
Mordell-Weil lattice is isomorphic to the {\em lattice} $\Hom
(C_1,C_2)$ with norm $\varphi \mapsto \deg (\varphi)$ up to torsion
(see \cite[Prop.3.1]{Shioda:CorrMWL}).
\end{example}

The next is the motivating example for studying the elliptic parameters and  the problems posed in \S1.1 in general.

\begin{example}[Inose's pencils] 
Using the twenty-four basic curves, we can find two disjoint divisors of Kodaira type
IV$^\ast$. Namely, take the following divisors shown in Fig.~\ref{J3}:
$$
\left\{
 \begin{array}{l} 
\Psi_1= G_1+G_2+G_3 +2(A_{01}+A_{02}+ A_{03}) + 3 F_0 ,\\
\Psi_2= F_1+F_2+F_3 +2(A_{10}+A_{20}+ A_{30}) + 3 G_0,
\end{array}
\right.
$$
\begin{figure}[!ht]
\begin{center}
\[
\includegraphics[scale=0.72]{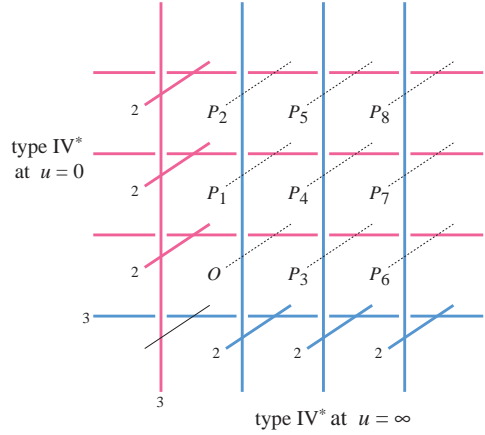}
\]
\caption{Inose's pencil (type $\J_{3}$)}
\label{J3}
\end{center}
\end{figure}
There is an elliptic fibration, called {\em Inose's pencil}, having
these divisors as the singular fibers over $ u=0$ and $u=\infty$, as
first shown by Inose \cite{Inose}. The elliptic parameter for this is
given by the function $u=t(=y_2/y_1) \in k(X)$, and the generic fiber
$E /k(t)$ is isomorphic to the cubic curve defined by the equation
$(\ref{Kummer})$ in the projective plane with inhomogeneous
coordinates $x_1,x_2$. (This belongs to type $\J_3$ in \cite{O}.)

It should be remarked that Kuwata \cite{Kuwata:MWrank} has succeeded
in constructing, by the use of Inose's pencil, some elliptic $K3$
surfaces with high Mordell-Weil rank which have an explicit defining
equation. For example, the base change $t=s^3$ gives rise to the
elliptic curve $E/k(s)$ which has the highest possible rank $r= 18$
(for $k=\C$) provided that $C_1$ and $C_2$ are mutually isogenous but
non-isomorphic elliptic curves with complex multiplications. We refer
to Kuwata \cite{Kuwata:MWrank} and Shioda \cite{K3SP},  
\cite{Shioda:CorrMWL} for more details including the defining equation
of $E$ in the Weiertrass form as well as the structure of MWL; see
also \S2.2.

\end{example}

\begin{example}
Besides the Kummer pencils (Example~1.1), the elliptic pencil on the
Kummer surface $X=\Km(C_{1}\times C_{2})$ which has been studied first
is perhaps the one introduced in Shioda-Inose \cite{IS}. It has
II$^{*}$, I$_{0}^{*}$, I$_{0}^{*}$ as reducible singular fibers (for
general values of $\l_{1}$ and $\l_{2}$). This has type $\J_{9}$ in
\cite{O} (see Fig.~\ref{J9}). Via the base change of degree~$2$, it
gives rise to an elliptic $K3$ surface with two II$^{*}$ fibers, which
plays an important role in the theory of singular $K3$ surfaces
\cite{IS} and which has been reconsidered by Morrison \cite{Morrison}
in a more general situation. It turns out that the elliptic parameter
and the defining equation for this type $\J_{9}$ is the hardest case
treated in this paper (see \S5.3).
\end{example}

\renewcommand{\baselinestretch}{1}
\begin{sidewaystable}
\footnotesize
\centering
\renewcommand\arraystretch{2}
\begin{tabular}{|c|c|c|c|}
\hline
Type 
&
\vtop {\hsize 3cm
Singular fibers
\vspace{6pt}}
& 
MWL 
& 
Elliptic parameter $u$ 
\\[5pt]\hline
$\J_{1}$ 
& 
$2\text{I}_{8} + 8\text{I}_{1} $ 
& 
$\Z^{2}\oplus \Z/2\Z$
& 
$\ds\frac{tx_{1}}{x_{2}}$ 
\\[5pt]\hline
$\J_{2}$ 
& 
$\text{I}_{4} +\text{I}_{12} + 8\text{I}_{1} $ 
&  
$A_2^*[2]\oplus \Z/2\Z$
& 
$\ds\frac{t(x_{1}-\l_{1})(x_{1}-x_{2})}{x_{2}(x_{2}-1)}$
\\[5pt]\hline
$\J_{3}$ 
& 
$2\text{IV}^{*} + 8\text{I}_{1}$
&  
$(A_2^*[2])^{2}$
& 
$t$
\\[5pt]\hline
$\J_{4}$ 
& 
$4\text{I}_{0}^{*}$ 
&  
$(\Z/2\Z)^{2}$
& 
$x_{i}$
\\[5pt]\hline
$\J_{5}$ 
& 
$\text{I}_{6}^{*} + 6\text{I}_{2}$ 
&  
$(\Z/2\Z)^{2}$
& 
$\ds\frac{(x_{1}-x_{2})\bigl(\l_{2}(x_{1}-\l_{1})+(\l_{1}-1)x_{2}\bigr)}
{(\l_{2}x_{1}-x_{2})\bigl(x_{1}-\l_{1}+(\l_{1}-1)x_{2}\bigr)}$
\\[5pt]\hline
$\J_{6}$ 
& 
2 $\text{I}_{2}^{*} + 4\text{I}_{2}$ 
&  
$(\Z/2\Z)^{2}$
& 
$\ds\frac{x_{1}}{x_{2}}$
\\[5pt]\hline
$\J_{7}$ 
& 
$\text{I}_{4}^{*} + 2\text{I}_{0}^{*} + 2\text{I}_{1}$ 
&  
$\Z/2\Z$
& 
$\ds\frac{(x_{2}-\l_{2})(x_{1}-x_{2})}{(x_{2}-1)(\l_{2}x_{1}-x_{2})}$
\\[5pt]\hline
$\J_{8}$ 
& 
$\text{III}^{*} + \text{I}_{2}^{*} + 3\text{I}_{2}+ \text{I}_{1}$
&  
$\Z/2\Z$
& 
$\ds-\frac{(x_{2}-\l_{2})(x_{1}-x_{2})} {\l_{2}(\l_{2}-1) x_{1}(x_{1}-1)}$
\\[5pt]\hline
$\J_{9}$
& 
$\text{II}^{*} + 2\text{I}_{0}^{*} + 2\text{I}_{1}$
&  
\{0\}
&
$\ds \frac{(x_{2}-\l_{2})(x_{1}-x_{2}) 
\bigl(\l_{2}x_{1} ( x_{1}-1)+(\l_{1}-1)(x_{2}-1)(\l_{2}x_{1}-x_{2})\bigr)}
{(x_{2}-1)(\l_{2}x_{1} -x_{2}) 
\bigl(\l_{2}x_{1} ( x_{1}-1 )+(\l_{1}-1)( x_{2}-\l_{2} )
(x_{1}-x_{2})\bigr)}$
\\[5pt]\hline
$\J_{10}$ 
&
$\text{I}_{8}^{*} + \text{I}_{0}^{*} + 4\text{I}_{1}$ 
&  
\{0\}
& $\ds\frac{( x_{2}-\l_{2} )(x_{1}-x_{2})
\bigl((\l_{1}-1)( x_{2}-1 )(\l_{2}x_{1} -x_{2}) +\l_{2}x_{1} ( x_{1}-1)\bigr)}
{x_{2}(x_{2}-1)(x_{1}-1)(\l_{2}x_{1} -x_{2})}$
\\[5pt]\hline
$\J_{11}$ 
& 
$2\text{I}_{4}^{*} + 4\text{I}_{1}$
&  
\{0\}
& 
$\ds\frac{x_{2}(x_{2}-\l_{2})(x_{1}-x_{2})}{x_{1}(x_{2}-1)(\l_{2}x_{1}-x_{2})}$
\\[5pt]\hline
\end{tabular}
\bigskip
\label{results}
\caption{Results}
\end{sidewaystable}

\subsection{Results}

In the following Table~1, we give a summary of the elliptic parameters
and the structure of the MWL for each type $\J_n$, to be
constructed in the subsequent sections.

The first column shows the type $\J_{n}$ of elliptic fibration
following Oguiso's notation (cf.~\cite{O}). The second column shows
the configuration of singular fibers in the \emph{generic case}, which
means that $\l_{1}$ and $\l_{2}$ are algebraically independent
elements of $k$ over $\Q_{0}$, where $\Q_{0}$ is the prime field in
$k$. The third column shows the structure of MWL of the generic fiber
$E$ over $k(u)$, again in the generic case. The last column gives the
elliptic parameter which can be used for any $\l_{1},\l_{2}(\neq
0,1)$.

The explicit form of defining equations should be found in the text,
since it is not suitable to tabulate here. We note that each of these
defining equations has coefficients in $\Q_{0}(\l_{1},\l_{2})(u)$,
where $u$ is the elliptic
parameter.

We see from the table that the elliptic parameters for $\J_n$ for
$n=1,2, 3$ are of the form $u=t \varphi(x_1,x_2)$ with
$\varphi(x_1,x_2) \in k(x_1,x_2)$, while those for $\J_n$ for $n>3$
are contained in $k(x_1,x_2)$.

\subsection{Basic strategy of construction}

Theoretically, constructing an elliptic fibration on a $K3$ surface is
to find a divisor that has the same type as a singular fiber in the
Kodaira's list (cf. \cite{Kodaira} \cite{P-S-S}). In practice,
however, we need to find two divisors, one for the fiber at $u=0$, and
the other for the fiber at $u=\infty$, to write down an actual
elliptic parameter. This is where the difficulty is.

Once an elliptic parameter is found, we want to find a change of
variables that converts the defining equation to a Weierstrass form.
In most cases we encounter an equation of the form
$y^{2}=(\text{quartic polynomial})$. We then use a standard algorithm
to transform it to a Weierstrass form (see for example Cassels
\cite{Cassels}, or Connell \cite{Connell}).

Some elliptic fibrations have nontrivial Mordell-Weil group. To
determine the structure of Mordell-Weil lattice, we can use the method
in \cite{Shioda:MWL} since we understand very well the intersection
between the section and the components of singular fibers.
Alternatively, we can compute the height pairing using the algorithm
in \cite{Kuwata:can-height} once we establish the conversion to the
Weierstrass form. Note that \cite{Shioda:MWL} and
\cite{Kuwata:can-height} use different normalization of the height
pairing, and they differ by a multiple of~$2$. In this article we
adopt the normalization used in~\cite{Shioda:MWL}.

\subsection{Remark}

Fix a type $\J_{n}$ $(n=1,\dots,11)$. As noted in \S1.5, each of the
defining equation of $E/k(u)$ constructed in this paper has the
coefficients in $\Q_{0}(\l_{1},\l_{2})(u)$, where $u$ is the elliptic
parameter, and $\l_{1}$ (resp.~$\l_{2}$) is the Legendre parameter for
$C_{1}$ (resp.~$C_{2}$). Given $C_{i}$, there are in general six
choices of $\l_{i}$ (i.e., six different level $2$-structures on
$C_{i}$). We have verified that, by different choices of $\l_{1}$ or
$\l_{2}$, we obtain as many nonisomorphic $E$'s belonging to the same
type $\J_{n}$, as predicted by Oguiso's result (\cite[Table~B,
p.~652]{O}), and thus solved Problem~2 when $C_{1}$ and $C_{2}$ are
not isogenous. The proof for this will be omitted in this paper, but
we write down the results in a single special case where we take
$C_1: y_1^2=x_1^3 -1$ and $C_2: y_2^2=x_2^3 - x_2$ (see \S6).

This paper is organized as follows.

\setcounter{tocdepth}{1}
\tableofcontents

\section{Elliptic parameters for $\J_{1}$, $\J_{3}$, $\J_{4}$, 
and $\J_{6}$}

In this section we construct elliptic fibrations that have two
singular fibers consisting only of the twenty-four basic curves. 
We use combinations of the following divisors of typical functions 
(cf. Examples in~\S1.4):
\begin{align*}
&(x_{1}) = 2 F_{1} + A_{10} + A_{11} + A_{12} + A_{13} 
- (2 F_{0} + A_{00} + A_{01} + A_{02} + A_{03}), \\
&(x_{2}) = 2 G_{1} + A_{01} + A_{11} + A_{21} + A_{31} 
- (2 G_{0} + A_{00} + A_{10} + A_{20} + A_{30}),\\
&(t) = G_{1} + G_{2} + G_{3} + 2 (A_{01} + A_{02} + A_{04}) + 3 F_{0}
\\
&\hspace{10em} - (F_{1} + F_{2} + F_{3}
+ 2 (A_{10} + A_{20} +  A_{30})+ 3 G_{0}).
\end{align*}

\subsection{$\J_{1}$}
\begin{figure}[!ht]
\begin{center}
\[
\includegraphics[scale=0.72]{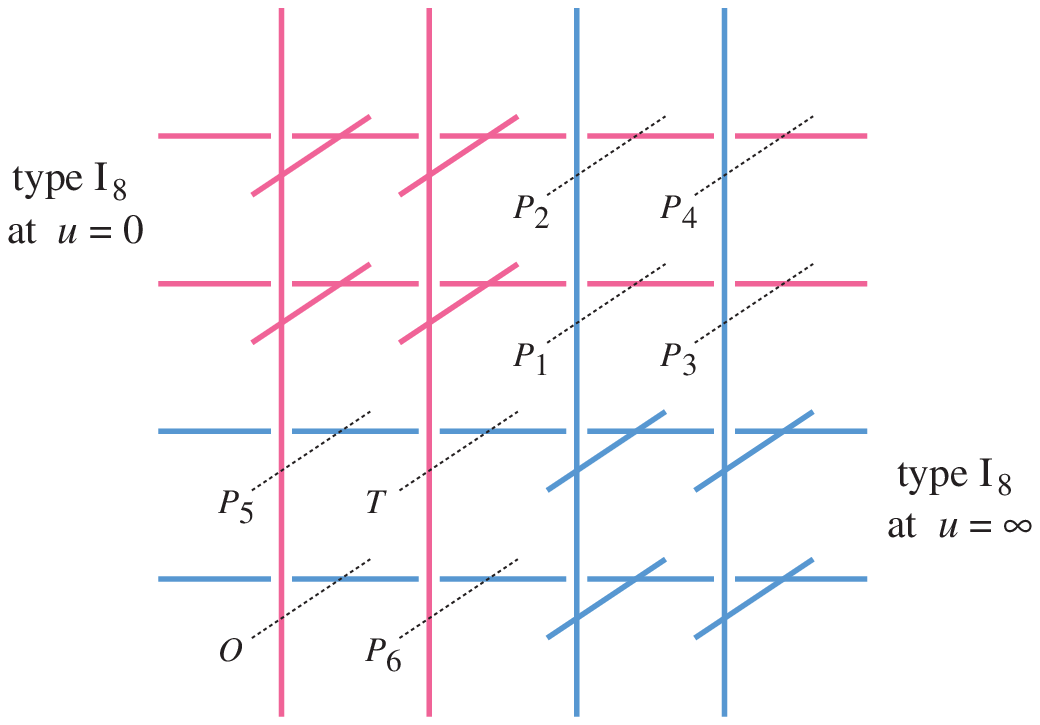}
\]
\caption{$\J_{1}$}
\label{J1}
\end{center}
\end{figure}

An elliptic parameter for the type $\J_{1}$ fibration is
given by
\[
u =\frac{tx_{1}}{x_{2}}.
\]
It is easy to verify that the divisor of $u$ is given by
\begin{multline*}
(u) = F_{0} + F_{1} + G_{2} + G_{3} + A_{02} + A_{03} + A_{12} + A_{13}
\\
- (G_{0} + G_{1} + F_{2} + F_{3} + A_{20} +A_{21} + A_{30} + A_{31}),
\end{multline*}
which is indicated in Fig.~\ref{J1}.
Choosing $A_{00}$ as the $0$-section of the group structure, we obtain
the Weierstrass equation of the elliptic fibration
\[
Y^2 = X^3 +\bigl((\l_{1}-1)^2
u^4-2(\l_{1}+1)(\l_{2}+1)u^2+(\l_{2}-1)^2\bigr) X^{2} 
+ 16 \l_{1}\l_{2}u^4X,
\]
where the change of variables is given by
\[
X =\frac{ 4t^2x_{1}^{3}}{x_{2}}, \quad Y = \frac{4t^2x_{1}^{3}
\bigl(t^{2 }x_{1}(x_{1}^{2}-\l_{1}) +x_{2}(x_{2}^2-\l_{2}) \bigr)}
{x_{2}^3}.
\]
Its discriminant is given by
\[
\Delta(u)=2^{12}\l_{1}^{2}\l_{2}^{2}\ u^{8}\,d(u)d(-u),
\]
where $d(u)$ is a polynomial of degree~$4$ in~$u$:
\begin{multline*}
d(u) = (\l_{1}-1)^{2}\,u^{4}+4(\l_{1}-1)\,u^{3} 
\\
-2(\l_{1}\l_{2}+\l_{1}+\l_{2}-3)\,u^{2}
+4(\l_{2}-1)\,u+ (\l_{2}-1)^{2}.
\end{multline*}
[The discriminant of $d(u)$ vanishes if and only if $\l_{1}=\l_{2}$.
If $\l_{1}=\l_{2}$, the elliptic fibration has two I$_{2}$ fibers for
general $\l_{1}$.]

The curve $A_{11}$ corresponds to the $2$-torsion section $T=(0,0)$.
The correspondence between the curves and the sections are as follows:
\[\renewcommand\arraystretch{1.41421356}
\begin{array}{ccl}
A_{22} & \leftrightarrow & P_{1} = \bigl(4u^2,
-4u^2((\l_{1}-1)u^{2}+\l_{2}-1)\bigr)
\\
A_{23} & \leftrightarrow & P_{2} = \bigl(4\l_{2}u^2,
-4\l_{2}u^2((\l_{1}-1)u^{2}-\l_{2}+1)\bigr)
\\
A_{32} & \leftrightarrow & P_{3} = \bigl(4\l_{1}u^2,
4\l_{1}u^2((\l_{1}-1)u^{2}-\l_{2}+1)\bigr)
\\
A_{33} & \leftrightarrow & P_{4} = \bigl(4\l_{1}\l_{2}u^2,
4\l_{1}\l_{2}u^2((\l_{1}-1)u^{2}+\l_{2}-1)\bigr)
\\
A_{01} & \leftrightarrow & P_{5} = \bigl(4\l_{2}, 4\l_{2}
((\l_{1}+1)u^{2}-\l_{2}-1)\bigr)
\\
A_{10} & \leftrightarrow & P_{6} = \bigl(4\l_{1}u^4,
-4\l_{1}u^4((\l_{1}+1)u^{2}-\l_{2}-1)\bigr)
\end{array}
\]
These sections satisfy the following relations.
\begin{gather*}
P_{3} = P_{2} + T,\qquad P_{4} = P_{1} + T,\\
P_{5} = P_{1} + P_{2},\qquad P_{6} = P_{5} + T.
\end{gather*}
The Mordell-Weil group is generated by $T$, $P_{1}$ and $P_{2}$ in the
general case where $C_{1}$ and $C_{2}$ are not isogenous. The height
matrix with respect to $\{P_{1}, P_{2} \}$ is shown to be
\[
\begin{pmatrix} 1 & 0 \\ 0 & 1 \end{pmatrix}.
\]

\subsection{$\J_{3}$}
As we have seen in~Example~1.2 (\S1.4),
\[
u = t
\]
gives an elliptic parameter of type $\J_{3}$.  
We regard \eqref{Kummer} as a cubic curve in $x_{1}$ and $x_{2}$ with
coefficients in $k(u)=k(t)$. We choose $(x_{1},x_{2})=(0,0)$ as the
origin of the group structure. The Weierstrass form is given by
\begin{align*}
Y^{2}&=X^{3} +4(\l_{1}+1)(\l_{2}+1)u^{2}X^{2} \\
&\quad+16u^{4}\bigl((\l_{1}\l_{2}+1)(\l_{1}+\l_{2}+1)-1\bigr)X \\
&\qquad
+16u^{4}\bigl((\l_{1}(\l_{1}-1)u^2+\l_{2}(\l_{2}-1))^{2}+4\l_{1}\l_{2}(\l_{1}+\l_{2})u^{2}
\bigr).
\end{align*}
(This is relatively simple, but the intermediate calculations are
rather complicated.) The change of variables between two forms of
equations is given by
\begin{align*}
&X = \frac{4
\bigl(\l_{2}(x_{1}-1)(x_{1}-\l_{1})+\l_{1}(x_{2}-1)
(x_{2}-\l_{2})-\l_{1}\l_{2}\bigr)\,t^{2}
}{x_{1} x_{2}},
\\
&Y = \frac{8(x_{2}-1)(x_{2}-\l_{2}) \bigl(\l_{2}(\l_{1}+1)x_{1} +
\l_{1}(\l_{2}+1)x_{2} - \l_{1}\l_{2} \bigr)\,t^{2}}{x_{1}^{2} x_{2}}\\
&\quad\qquad
+\frac{4\l_{1}\bigl((\l_{1}+1)x_{1}-2\l_{1}\bigr)\,t^{4}}{x_{1}}
+\frac{4\l_{2}\bigl((\l_{2}+1)x_{2}-2\l_{2}\bigr)\,t^{2}}{x_{2}}.
\end{align*}
The discriminant is of the form
\(
u^{8} d(u),
\)
where $d(u)$ is an irreducible polynomial of degree~$8$. Besides two
IV$^{*}$ fibers, the elliptic fibration has eight I$_{1}$ fibers in
the generic case. These eight I$_{1}$ fibers can degenerate in four
different ways; 2 I$_{2}$ + 4 I$_{1}$, 4 I$_{2}$, 4 II or 2 IV. For
more detail, see Prop.~5.1 in \cite{Shioda:CorrMWL}.

There are eight other $A_{ij}$'s which define sections;
the correspondence between these curves and the sections is as follows:
\[\renewcommand\arraystretch{1.41421356}
\setlength{\arraycolsep}{3pt}
\begin{array}{ccl}
A_{12} & \leftrightarrow & P_{1} =
\bigl(4u^{2}(\l_{1}^{2}u^{2}-\l_{2}(\l_{1}+1)), \\
&& \qquad\qquad
-4u^{2}(2\l_{1}^{3}u^{4}-\l_{1}(\l_{1}+1)(2\l_{2}-1)u^{2}
+\l_{2}(\l_{2}-1))\bigr)
\\[4pt]
A_{13} & \leftrightarrow & P_{2} =
\left(\dfrac{4u^{2}(\l_{1}^{2}u^{2}-\l_{2}^{2}(\l_{1}+1))}{\l_{2}^{2}},
\right.\\[8pt]
&& \qquad\quad\left.
-\dfrac{4u^{2}(2\l_{1}^{3}u^{4}+\l_{1}\l_{2}^{2}
(\l_{1}+1)(\l_{2}-2)u^{2}-\l_{2}^{4}(\l_{2}-1))}
{\l_{2}^{3}}\right)
\\
A_{21} & \leftrightarrow & P_{3} = \bigl(-4(\l_{1}(\l_{2}+1)u^{2}-\l_{2}^{2}), \\
&& \qquad\qquad -4(\l_{1}(\l_{1}-1)u^{4}
-\l_{2}(\l_{2}+1)(2\l_{1}-1)u^{2}+2\l_{2}^{3})\bigr)
\\
A_{22} & \leftrightarrow & P_{4} = \bigl(-4\l_{1}\l_{2}u^{2},
-4u^{2}(\l_{1}(\l_{1}-1)u^{2}+\l_{2}(\l_{2}-1))\bigr)
\\
A_{23} & \leftrightarrow & P_{5} = \bigl(-4\l_{1}u^{2},
-4u^{2}(\l_{1}(\l_{1}-1)u^{2}-\l_{2}(\l_{2}-1))\bigr)
\\[4pt]
A_{31} & \leftrightarrow & P_{6} =
\left(-\dfrac{4(\l_{1}^{2}(\l_{2}+1)u^{2}-\l_{2}^{2})}{\l_{1}^{2}},
\right.\\[8pt]
&& \qquad\qquad \left.
\dfrac{4(\l_{1}^{4}(\l_{1}-1)u^{4}-\l_{1}^{2}\l_{2}(\l_{1}-2)(\l_{2}+1)u^{2}-2\l_{2}^{3})}
{\l_{1}^{3}}\right)
\\
A_{32} & \leftrightarrow & P_{7} = \bigl(-4\l_{2}u^{2},
4u^{2}(\l_{1}(\l_{1}-1)u^{2}-\l_{2}(\l_{2}-1))\bigr)
\\
A_{33} & \leftrightarrow & P_{8} = \bigl(-4u^{2},
4u^{2}(\l_{1}(\l_{1}-1)u^{2}+\l_{2}(\l_{2}-1))\bigr)
\end{array}
\]
These sections satisfy the following relations:
\begin{align*}
P_{1} = P_{5} + P_{8}, \qquad P_{2} = P_{4} + P_{7}, \\
P_{3} = P_{7} + P_{8},\qquad P_{6} = P_{4} + P_{5}.
\end{align*}
We can show that $P_{4}, P_{8},P_{5}$, and $P_{7}$ generate the Mordell-Weil
group in the generic case. The height matrix with respect to the basis
$\{P_{4}, P_{8},P_{5}, P_{7}\}$ is
\[\renewcommand{\arraystretch}{1.2}\setlength{\arraycolsep}{6pt}
\begin{pmatrix} \frac{4}{3}& \frac{2}{3}  & 0 & 0 \\[3pt]
\frac{2}{3} & \frac{4}{3} & 0 & 0\\[3pt]
0 & 0 & \frac{4}{3} & \frac{2}{3}\\[3pt]
0 & 0 &\frac{2}{3} & \frac{4}{3}
\end{pmatrix}.
\]
This is the direct sum of two copies of $A_{2}^{*}[2]$, the dual
lattice of $A_{2}$ scaled by~$2$.

\subsection{$\J_{4}$}
The elliptic parameter for the fibration $\pi_{1}$ in~Example~1.1 is
given by
\[
u = x_{1},
\]
while the elliptic parameter for $\pi_{2}$ is given by
$u=x_{2}$.
For $\pi_{1}$, the change of variables
\begin{align*}
&X = u(u-1)(u-\l_{1})x_{2},
\\
&Y = u^{2}(u-1)^{2}(u-\l_{1})^{2}\,t,
\end{align*}
converts the equation \eqref{Kummer} to
\[
Y^{2}=X\bigl(X-u(u-1)(u-\l_{1})\bigr)\bigl(X-\l_{2}\,u(u-1)(u-\l_{1})\bigr).
\]
The curve $G_{0}$ is the $0$-section.  Other sections are:
\[\renewcommand{\arraystretch}{1.41421356}
\begin{array}{ccl}
G_{1} & \leftrightarrow & T_{1} = (0,0),
\\
G_{2} & \leftrightarrow & T_{2} = (u(u-1)(u-\l_{1}), 0),
\\
G_{3} & \leftrightarrow & T_{3} = (\l_{2}\,u(u-1)(u-\l_{1}), 0).
\end{array}
\]
Similar results hold for $\pi_{2}$.

\subsection{$\J_{6}$}
The divisor of the function $x_{1}/x_{2}$ is given by
\begin{multline*}
\Bigl(\frac{x_{1}}{x_{2}}\Bigr) 
= 2 (F_{1} + A_{10}  + G_{0}) + A_{12} + A_{13} + A_{20} + A_{30}\\
- \bigl(2 (F_{0}  +  A_{01} + G_{1}) + A_{02} + A_{03} + A_{21} + A_{31}\bigr).
\end{multline*}
This is the difference of two disjoint divisors of type I$_{2}^{*}$, and thus 
\[
u =\frac{x_{1}}{x_{2}}
\]
is an elliptic parameter of type $\J_{6}$. 

\begin{figure}[!ht]
\begin{center}
\[
\includegraphics[scale=0.72]{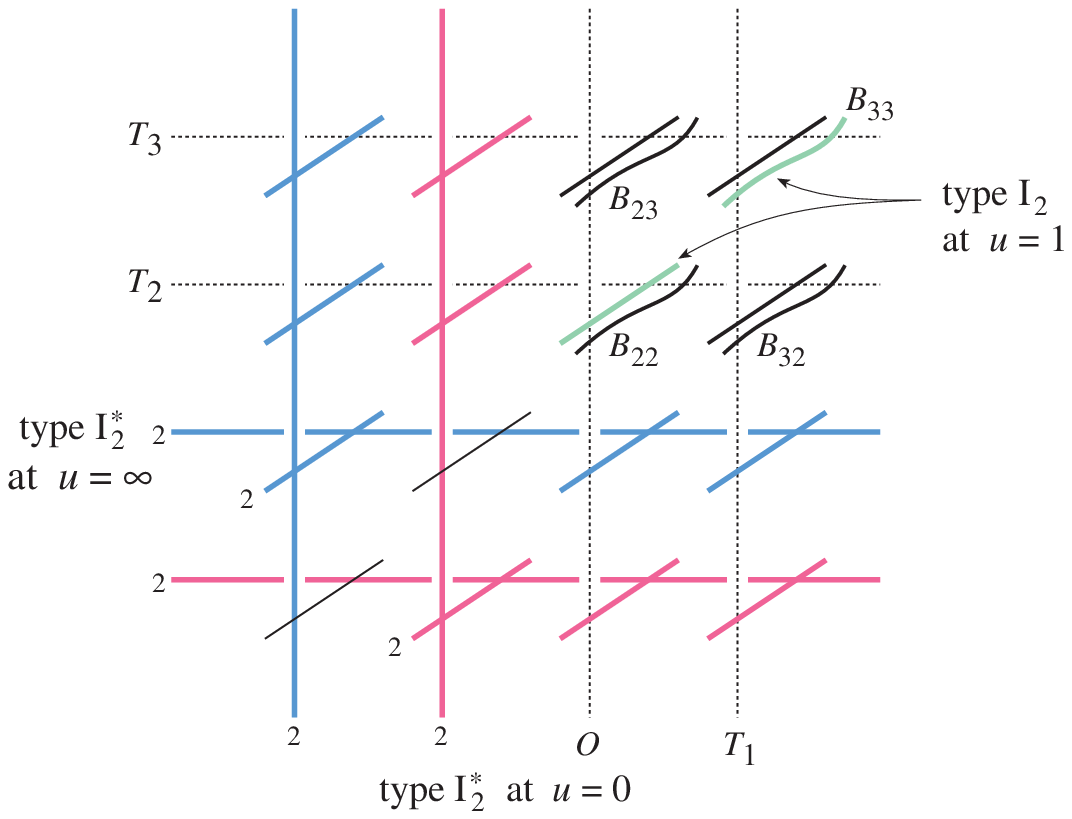}
\]
\caption{$\J_{6}$}
\label{J6}
\end{center}
\end{figure}

In order to write down a Weierstrass equation using the curve $F_{2}$
as the $0$-section, we put
\begin{align*}
X &=\frac{x_{1} (x_{1}-\l_{1})
(x_{1}-x_{2})(\l_{2}x_{1}-x_{2})}{(x_{1}-1)x_{2}^3},
\\
Y &= \frac{(\l_{1}-1)\, t\, x_{1}^3(x_{1}-\l_{1}) (x_{1}-x_{2}) (\l_{2}
x_{1}-x_{2})}{(x_{1}-1)x_{2}^5}.
\end{align*}
Then we obtain the Weierstrass equation
\begin{equation*}
Y^2 = X\bigl(X-u(u-1)(\l_{2}u-\l_{1})\bigr)
\bigl(X-u(u-\l_{1})(\l_{2}u-1)\bigr).
\end{equation*}
Its discriminant is given by
\[
\Delta(u) = 16{u}^{8}(\lambda_{1}-1)^{2} (\lambda_{2}-1)^{2} 
(u-1)^{2} (u-\lambda_{1}) ^{2} (\lambda_{2}u-1)^{2} 
(\lambda_{2}u-\lambda_{1})^{2}.
\]
Besides two I$_{2}^{*}$ fibers at $u=0$ and $\infty$, there are four
I$_{2}$ fibers at $u=1,\l_{1}, 1/\l_{2}$ and $\l_{1}/\l_{2}$. This
elliptic surface has the following three $2$-torsion sections:
\[\renewcommand{\arraystretch}{1.41421356}
\begin{array}{ccl}
F_{3} & \leftrightarrow & T_{1} = (0,0),
\\
G_{2} & \leftrightarrow & T_{2} = \bigl(u(u-\l_{1})(\l_{2}u-1),0\bigr),
\\
G_{3} & \leftrightarrow & T_{3} = \bigl(u(u-1)(\l_{2}u-\l_{1}), 0\bigr),
\end{array}
\]

Note that $A_{22}$, $A_{23}$, $A_{32}$, and $A_{33}$ are components of
four I$_{2}$ fibers. The other components of these four I$_{2}$ fibers
are new $(-2)$-curves not among the basic curves, which will be
clarified in \S3.2.

\section{More $(-2)$-curves}

In order to describe elliptic parameters for other types, we need more
$(-2)$-curves than the basic curves. When we constructed elliptic
parameters of type $\J_{6}$ just above, we obtained some new
$(-2)$-curves as components of I$_{2}$ fibers. In this section we give
a systematic way to obtain such $(-2)$-curves.

For our purpose, it is convenient to regard $X=\Km(C_{1}\times C_{2})$
as a double cover of the product of projective lines: $\Proj^{1}\times
\Proj^{1}=\{(x_{1}:z_{1}),(x_{2}:z_{2})\}$. Let $p_{i}:C_{i}\to
\Proj^{1}$ ($i=1,2$) be the projection given by
\[
\begin{array}{rccc}
p_{i}:&C_{i} & \longrightarrow & \Proj^{1}\\[3pt]
&(x_{i}:y_{i}:z_{i}) & \longmapsto & 
\begin{cases} (x_{i}:z_{i}) & \text{if $z_{i}\neq 0$}
\\[3pt] 
(1:0) &\text{if $z_{i}=0$} \end{cases}
\end{array}
\]
Then the map $p_{1}\times p_{2} : A=C_{1}\times C_{2} \to
\Proj^{1}\times\Proj^{1}$ factors through $\bar\pi:A/\iota_{A} \to
\Proj^{1}\times\Proj^{1}$. Let $\pi$ be the morphism of degree two
from $X$ to $\Proj^{1}\times\Proj^{1}$ that makes the following
diagram commutative:
\[
\begin{array}{rcccc}
&& X && \\
&&\downarrow&\overset{\hbox{\scriptsize $\pi$}}\searrow& \\[4pt]
A &\longrightarrow & A/\iota_{A}
&\underset{\bar\pi}{\longrightarrow} &\Proj^{1}\times\Proj^{1}
\end{array}
\]
We denote by $R_{ij}$ the point in $\Proj^{1}\times \Proj^{1}$ that is
the image of the exceptional curve $A_{ij}$ by $\pi$. To obtain more
$(-2)$-curves, we look for curves in $\Proj^{1}\times \Proj^{1}$ which
lift to a $(-2)$-curve via the map $\pi$.

\subsection{$(1,1)$-curves}
Let $L$ be a curve in $\Proj^{1}\times \Proj^{1}$ defined by a
bihomogeneous equation of bidegree $(1,1)$:
\[
a x_{1}x_{2} + b x_{1}z_{2} + c z_{1}x_{2} + d z_{1}z_{2} =0.
\]
We call such a curve $(1,1)$-curve for short. By an abuse of
notation, we denote the image of $F_{i}$ and $G_{i}$ under $\pi:S\to
\Proj^{1}\times \Proj^{1}$ by the same letters $F_{i}$ and $G_{i}$,
respectively.   For example, $F_{1}$ is the curve with the equation 
$x_{1}=0$, and $G_{3}$ with $x_{2}-\l_{2}z_{2}=0$, etc.

Let $L$ be a $(1,1)$-curve in $\Proj^{1}\times \Proj^{1}$. Its
pullback $\pi^{-1}(L)$ ramifies at the intersection of $L$ and $F_{i}$
or $G_{j}$, except when the intersection point falls on
$R_{ij}=F_{i}\cap G_{j}$.

\begin{lem}
Let $L$ be a $(1,1)$-curve. Then, 
\begin{enumerate}
\item If $L$ passes three of sixteen $R_{ij}$'s, then $\pi^{-1}(L)$ is 
a curve of genus~$0$.
\item If $L$ passes two out of sixteen $R_{ij}$'s, then $\pi^{-1}(L)$ 
is a curve of genus~$1$.
\end{enumerate}
\end{lem}

\begin{proof}
In general a $(1,1)$-curve $L$ intersects with $\sum F_{i}$ (resp.
$\sum G_{j}$) at four points. If $L$ passes three of sixteen
$R_{ij}$'s, then it intersects with $F_{i}$ one more time and $G_{j}$
one more time outside $R_{ij}$. This implies that $\pi^{-1}(L)$
ramifies at two points. By Hurwitz's theorem $\pi^{-1}(L)$ is a curve
of genus~$0$. Similarly, if $L$ passes two out of sixteen $R_{ij}$'s,
$\pi^{-1}(L)$ ramifies at four points, and it is a curve of genus~$1$.
\end{proof}
 
A $(1,1)$-curve is uniquely determined by a set of three points in a
general position. If we choose $R_{i_{0}j_{0}}$,
$R_{i_{1}j_{1}}$, $R_{i_{2}j_{2}}$ so that no two of them are on the
same $F_{i}$ or $G_{j}$, then they are in general position. Let
$i_{3}$ and $j_{3}$ be the missing indices. In other words, we choose
$i_{3}$ and $j_{3}$ such that
$\{i_{0},i_{1},i_{2},i_{3}\}=\{j_{0},j_{1},j_{2},j_{3}\}=\{0,1,2,3\}$.
Under the condition that the two elliptic curves $C_{1}$ and $C_{2}$
are not isomorphic, the $(1,1)$-curve passing through
$R_{i_{0}j_{0}}$, $R_{i_{1}j_{1}}$, and $R_{i_{2}j_{2}}$ does not pass
$R_{i_{3}j_{3}}$. Thus, choosing $R_{i_{0}j_{0}}$, $R_{i_{1}j_{1}}$,
$R_{i_{2}j_{2}}$ we obtain a $(1,1)$-curve whose pullback by $\pi$ is
an irreducible $(-2)$-curve in $X$. We denote such a $(1,1)$-curve by
$L_{i_{0}j_{0},i_{1}j_{1},i_{2}j_{2}}$, and its pullback by $\tilde
L_{i_{0}j_{0},i_{1}j_{1},i_{2}j_{2}}$. There are ninety-six such
$(-2)$-curves. Also note that $\tilde
L_{i_{0}j_{0},i_{1}j_{1},i_{2}j_{2}}$ intersects twice with each of
$A_{i_{0}j_{0}}$, $A_{i_{1}j_{1}}$, and $A_{i_{2}j_{2}}$.

The $(1,1)$-curve $L_{00,11,22}$ passes through $R_{00}$, $R_{11}$,
$R_{22}$. It is given by the bihomogeneous equation
$x_{2}z_{1}-x_{1}z_{2}=0$. In the sequel we write it in the affine
form $x_{2}-x_{1}=0$ for simplicity. $\tilde L_{00,11,22}$ is denoted
by $A^{44}$ in Oguiso \cite{O}, which appears in the $\J_{2}$
fibration. We denote it by $B_{33}$ to make it consistent with our
notation, indicating that it intersects with $F_{3}$ and $G_{3}$
outside $A_{33}$. Note, however, that there are six $(-2)$-curves of
the form $\tilde L_{i_{0}j_{0},i_{1}j_{1},i_{2}j_{2}}$ that intersect
with $F_{3}$ and $G_{3}$.

\begin{figure}[!ht]
\begin{center}
\[
\includegraphics[scale=0.8]{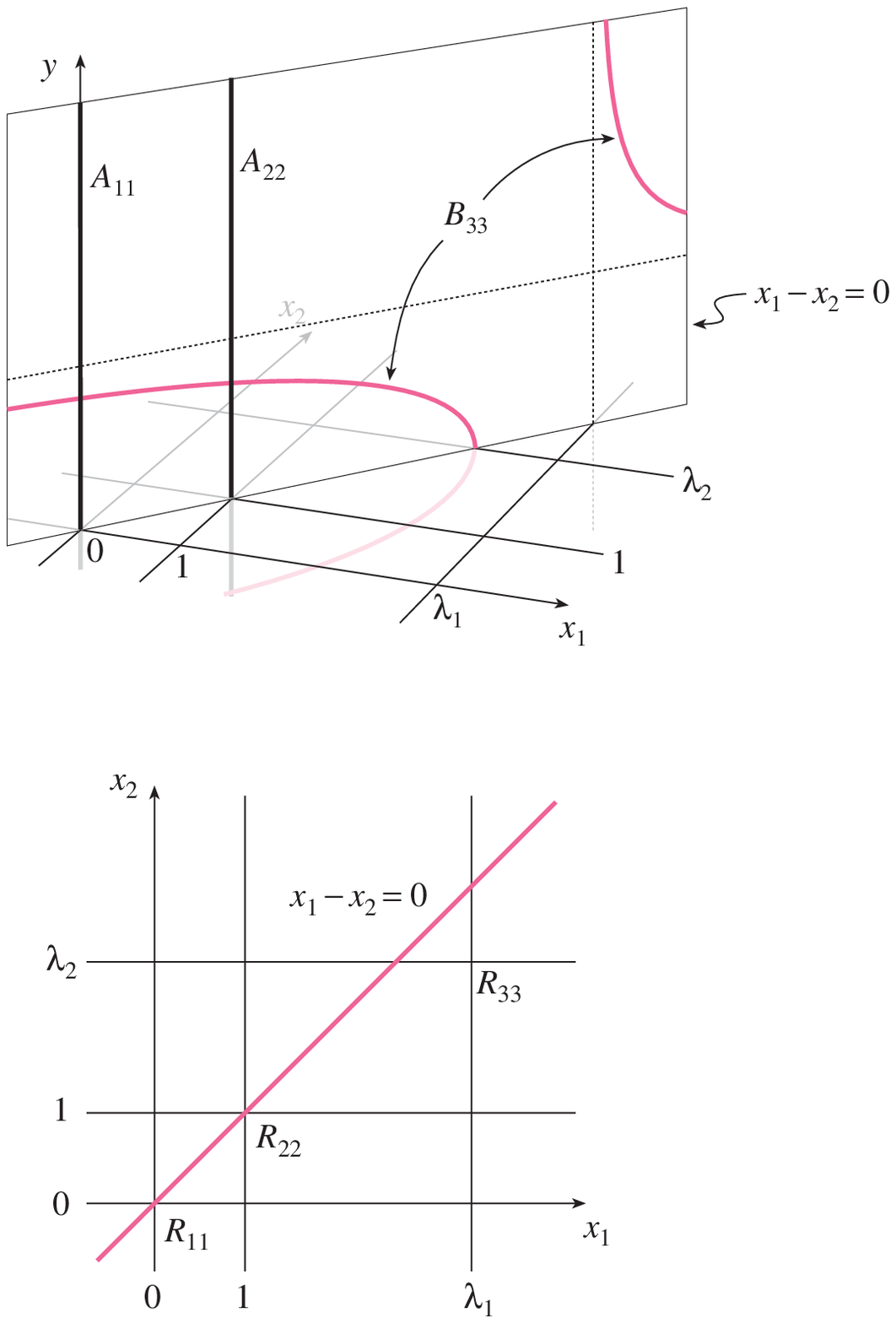}
\]
\caption{$(-2)$-curve $B_{33}$}
\label{B_33}
\end{center}
\end{figure}

Fig.~\ref{B_33} shows the curve $B_{33}$ in the affine space
$\A_{x_{1}}\times\A_{x_{2}}\times \A_{t}$. As a matter of fact, if we
substitute $x_{2}$ by $x_{1}$ in \eqref{Kummer}, the equation
factorizes into
 \[
x_{1} (x_{1}-1) (x_{1} t^2-x_{1}-t^2 \l_{1}+ \l_{2}) = 0,
\]
which implies that the intersection between $x_{2}-x_{1}=0$ and the
affine Kummer surface \eqref{Kummer} has three irreducible components,
namely $A_{11}$, $A_{22}$, and $B_{33}$.  We also see that a
parametrization of $B_{33}$ is given by
\[
(x_{1},x_{2},t)=
\left(\frac{\l_{1}s^{2}-1}{s^{2}-1},\frac{\l_{1}s^{2}-1}{s^{2}-1},s\right).
\]

The zero divisor of the function $x_{2}-x_{1}\in k(x_{1},x_{2},t)$ is
$A_{11} + A_{22} + B_{33}$, while the polar divisor is of the form
$D_{1} + D_{2} + r A_{00}$, where
\begin{align*}
D_{1} = 2 F_{0} +  A_{00} + A_{01} + A_{02} + A_{03}, \\
D_{2} = 2 G_{0} + A_{00}  + A_{10} +A_{20} + A_{30}.
\end{align*}
Since $A_{00}$ intersects twice with the divisor $A_{11} + A_{22} +
B_{33}$, the intersection number $A_{00}\cdot(D_{1} + D_{2} + r
A_{00}) $ must be~$2$, which implies $r=-1$.
This shows 
\begin{multline*}
(x_{2}-x_{1}) = A_{11} + A_{22} + B_{33} \\
- (2 F_{0} + 2 G_{0} + A_{00} + A_{01} + A_{02} + A_{03} + A_{10} +
A_{20} + A_{30}).
\end{multline*}
This and similar calculations of divisors are used to find the
elliptic parameter with a prescribed divisor in \S4 and \S5.

\subsection{I$_{2}$ fibers of type $\J_{6}$ fibration}
\begin{figure}[!hbt]
\begin{center}
\[
\includegraphics[scale=0.8]{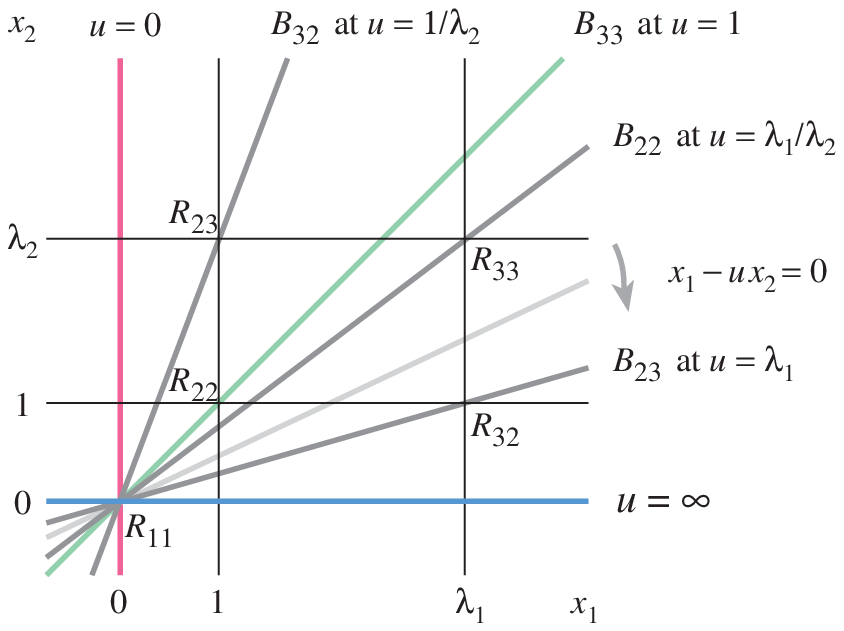}
\]
\caption{pencil of $(1,1)$-curves}
\label{J6_plane}
\end{center}
\end{figure}

The elliptic parameter $u=x_{1}/x_{2}$, which is of type $\J_{6}$,
defines a pencil of $(1,1)$-curves $x_{1}-ux_{2}=0$. The general fiber
of this elliptic fibration is the pullback of a $(1,1)$-curve passing
through $R_{00}$ and $R_{11}$. If $x_{1}-ux_{2}=0$ passes through a
third $R_{ij}$, then its pullback is a singular fiber (see
Fig.~\ref{J6_plane}).
Four fibers of type I$_{2}$, which are mentioned in \S2.4 arise as
follows:
\[
\begin{array}{llcll}
\tilde L_{00,11,23}+A_{23} &\text{ at } u=1/\l_{2},&&
B_{33}+A_{22}  & \text{ at }  u=1,  \\
\tilde L_{00,11,33}+A_{33} &\text{ at } u=\l_{1}/\l_{2}, &&
\tilde L_{00,11,32}+A_{32} &\text{ at } u=\l_{1}.
\end{array}
\]

\subsection{Notation}\label{notation}
Even though the notation ``$B_{33}$'' is ambiguous as we mentioned
earlier, it is quite convenient. We thus use the following notation in
the sequel:
\begin{equation}
\begin{array}{ll}
   B_{32} =\tilde L_{00,11,23}: \l_{2}x_{1}-x_{2} = 0,&
   B_{33} =\tilde L_{00,11,22}: x_{1}-x_{2} = 0, \\ 
   B_{22} =\tilde L_{00,11,33}: \l_{2}x_{1}-\l_{1}x_{2} = 0, &
   B_{23} =\tilde L_{00,11,32}:x_{1}-\l_{1}x_{2} = 0, .
\end{array}
\end{equation}
Later in \S4.3 and \S5.2, we introduce more $(-2)$-curves of this
type, $B_{31}$, $B_{12}$ and $B_{13}$.

\section{Elliptic parameters for $\J_{2}$, $\J_{7}$, $\J_{8}$
 and $\J_{11}$}

\subsection{$\J_{2}$}
Using $B_{33}$, we can construct an elliptic parameter of type
$\J_{2}$. In fact, the divisor
\[
\Psi_{2,0}=F_{3}+A_{33}+G_{3}+B_{33}
\]
is a divisor of type I$_{4}$, and it does not intersect with the
divisor of type I$_{12}$ given by
\begin{multline*}
\Psi_{2,\infty}=F_{0} + A_{02} + G_{2} + A_{12} + F_{1} + A_{10} \\
+ G_{0} + A_{20} + F_{2} + A_{21} + G_{1} + A_{01}
\end{multline*}
(see Fig.~\ref{J2} below). It turns out that the divisor of the
function
\[
u =\frac{t(x_{1}-\l_{1})(x_{1}-x_{2})}{x_{2}(x_{2}-1)}
\]
is $\Psi_{2,0}-\Psi_{2,\infty}$, and it is an elliptic parameter of
type $\J_{2}$.
\begin{figure}[!ht]
\begin{center}
\[
\includegraphics[scale=0.72]{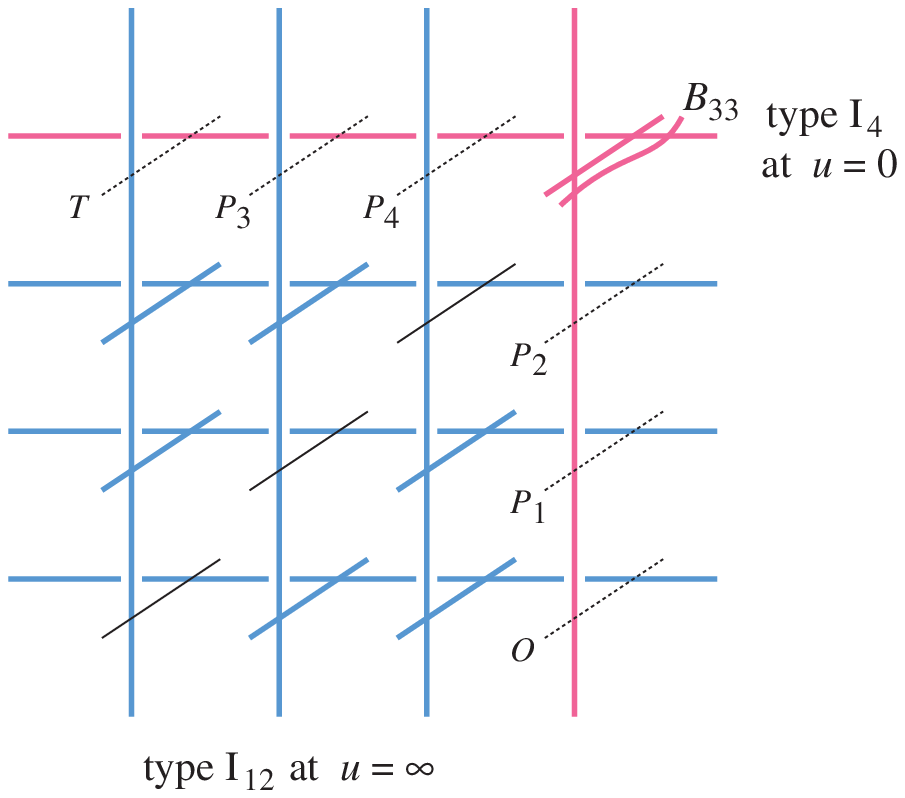}
\]
\caption{$\J_{2}$}
\label{J2}
\end{center}
\end{figure}
Choosing $A_{30}$ as the $0$-section, we obtain
the Weierstrass equation
\begin{multline*}
Y^2 = X^{3} + \bigl(
u^4+2\,(2\l_{1}\l_{2}-\l_{1}-\l_{2}+2)\,u^2+(\l_{2}-\l_{1})^2\bigr)X^{2} \\
 -16\l_{1}\l_{2}(\l_{1}-1)(\l_{2}-1)\,u^2X,
\end{multline*}
where the change of variables is given by
\begin{align*}
&X =
-\frac{4\l_{1}(\l_{1}-1)(x_{1}-x_{2})(x_{2}-\l_{2})}{x_{1}(x_{1}-1)},
\\
&Y =
-\frac{4\l_{1}(\l_{1}-1)(x_{1}-x_{2})(x_{2}-\l_{2})(2x_{1}-2x_{2}-\l_{1}+\l_{2})
} {x_{1}(x_{1}-1)}
\\
&\qquad\qquad +
\frac{4\l_{1}(\l_{1}-1)(x_{1}-x_{2})^{2}(x_{2}-\l_{2})^{3}
(2x_{1}x_{2}-x_{1}-x_{2})} {t^{2}x_{1}^{3}(x_{1}-1)^{3}}.
\end{align*}
The discriminant of the fibration is of the form $u^{4}d(u)$, where
$d(u)$ is a polynomial of degree~$8$. The discriminant of $d(u)$
vanishes if and only if
\[
\l_{2}=\l_{1}, 1-\l_{1}, \frac{1}{\l_{1}}, \text{or } \frac{\l_{1}}{\l_{1}-1}.
\]

The curve $A_{03}$ corresponds to the $2$-torsion section $T=(0,0)$.
The correspondence between the curves and the sections is as follows:
\[\renewcommand\arraystretch{1.41421356}
\begin{array}{ccl}
A_{31} & \leftrightarrow & P_{1} = \bigl(4\l_{1}\l_{2},
4\l_{1}\l_{2}(u^{2}+\l_{1}+\l_{2}))\bigr)
\\
A_{32} & \leftrightarrow & P_{2} = \bigl(4(\l_{1}-1)(\l_{2}-1),
\\ && \qquad\qquad
-4(\l_{1}-1)(\l_{2}-1)(u^{2}-\l_{1}-\l_{2}+2))\bigr)
\\
A_{13} & \leftrightarrow & P_{3} = \bigl(-4u^{2}(\l_{1}-1)(\l_{2}-1),
\\ && \qquad\qquad\quad
4(\l_{1}-1)(\l_{2}-1)u^{2}(u^{2}+\l_{1}+\l_{2}))\bigr)
\\
A_{23} & \leftrightarrow & P_{4} = \bigl(-4\l_{1}\l_{2}u^{2},
-4\l_{1}\l_{2}u^{2}(u^{2}-\l_{1}-\l_{2}+2))\bigr)
\end{array}
\]
These sections satisfy the following relations.
\[
P_{3} = P_{1} + T,\qquad P_{4} = P_{2} + T.
\]
The Mordell-Weil group is generated by $T$, $P_{1}$ and $P_{2}$ in the
general case where $C_{1}$ and $C_{2}$ are not isogenous. The height
matrix with respect to $\{P_{1}, P_{2} \}$ is shown to be
\[\renewcommand{\arraystretch}{1.2}
\begin{pmatrix} \frac{4}{3} & \frac{2}{3} \\[3pt] \frac{2}{3} &
\frac{4}{3} \end{pmatrix}.
\]
Thus the Mordell-Weil lattice is isomorphic to $A_{2}^{*}[2]$.

\subsection{$\J_{7}$} 
\begin{figure}[!ht]
\begin{center}
\[
\includegraphics[scale=0.72]{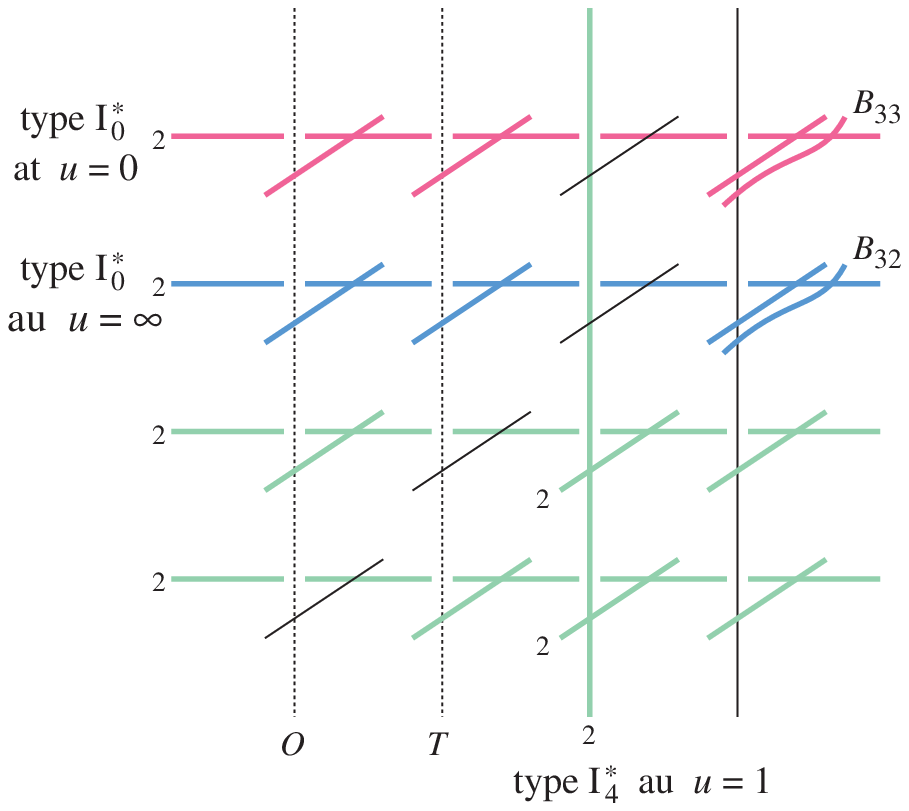}
\]
\caption{$\J_{7}$}
\label{J7}
\end{center}
\end{figure}

Using the curves $B_{33}$ and $B_{32}$ introduced in \S\ref{notation}, 
we can form two disjoint divisors of type I$_{0}^{*}$:
\begin{align*}
& \Psi_{7,0} = 2 G_{3} + A_{03} + A_{13} + A_{33} + B_{33}, \\
& \Psi_{7,\infty} =2 G_{2} + A_{02} + A_{12} + A_{32} + B_{32}.
\end{align*}
Looking for the function whose divisor is $\Psi_{7,0}
-\Psi_{7,\infty}$, we obtain the elliptic parameter
\begin{equation}\label{ell_param:J_7}
u =\frac{(x_{2}-\l_{2})(x_{1}-x_{2})}{(x_{2}-1)(\l_{2}x_{1}-x_{2})}.
\end{equation}
The divisor of the function
\[
u -1 = -\frac{(\l_{2}-1) \,x_{2}(x_{1}-1)}{(x_{2}-1)(\l_{2}x_{1}-x_{2})}
\]
is given by the following divisor consisting only of the basic curves:
\[
A_{01} + A_{31} + 2 G_{1} + 2 A_{21} + 2 F_{2} + 2 A_{20} + 2 G_{0} +  A_{10} +  A_{30}.
\]
This is a singular fiber of type I$_{4}^{*}$ (see Fig.~\ref{J7}).
Thus, the elliptic parameter given by~\eqref{ell_param:J_7} is of type
$\J_{7}$.

The change of variables
\begin{align*}
X &=\frac{\l_{2}u(u-1)^2x_{1}}{x_{2}} \\
&=\frac{\l_{2}(\l_{2}-1)^2 x_{1} (x_{1}-1)^{2} x_{2} (x_{2}-\l_{2})
(x_{1}-x_{2})} {(x_{2}-1)^3(\l_{2}x_{1}-x_{2})^3},
\\
Y &=\frac{\l_{2}(\l_{2}-1)u^{2}(u-1)^2}{t}\\
&=\frac{\l_{2}(\l_{2}-1)^3 (x_{1}-1)^2 x_{2}^2 (x_{2}-\l_{2})^2
(x_{1}-x_{2})^2} {t(x_{2}-1)^{4} (\l_{2}x_{1}-x_{2})^{4}},
\end{align*}
converts \eqref{Kummer} to the Weierstrass equation
\[
Y^2 = X^3 -u(u-1)\bigl((\l_{1}\l_{2}+1)u
-\l_{1}-\l_{2}\bigr)X^2+\l_{1}\l_{2}u^2(u-1)^4X.
\]
Its discriminant is of the form $u^{6}(u-1)^{10}d(u)$, where $d(u)$ is
a polynomial of degree~$2$.

Generically, it has only one section other than $0$-section:
\[
F_{2}  \quad\leftrightarrow\quad  T = (0,0).
\]

\subsection{$\J_{8}$}
To find an elliptic parameter of type $\J_{8}$, we need to construct a
I$_{2}^{*}$ fiber. For this, we can make use of $B_{33}$ once again.
The divisor
\[
\Psi_{8,0} = A_{01} + A_{02} + 2 F_{0} + 2 A_{03} + 2 G_{3} + A_{33} + B_{33}
\]
is of type I$_{2}^{*}$ and it does not intersect with the divisor
\[
\Psi_{8,\infty} = A_{12} + 2 F_{1} + 3 A_{10} + 4 E_{0} + 3 A_{20} + 2 F_{2} + A_{21} + 2 A_{30},
\]
which is of type III$^{*}$. We look for a function whose divisor is
$\Psi_{8,0} - \Psi_{8,\infty}$, and we obtain the elliptic parameter
of type $\J_{8}$
\[
u =-\frac{(x_{2}-\l_{2})(x_{1}-x_{2})} {\l_{2}(\l_{2}-1) x_{1}(x_{1}-1)}.
\]
\begin{figure}[!ht]
\begin{center}
\[
\includegraphics[scale=0.72]{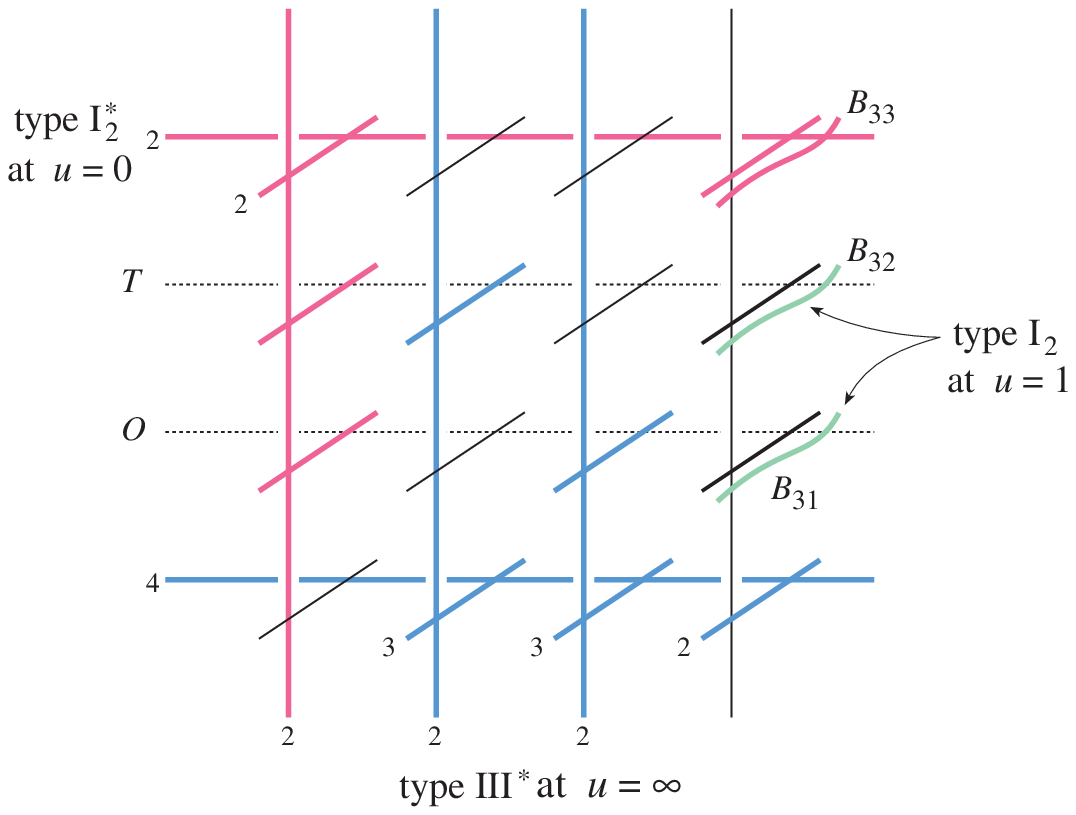}
\]
\caption{$\J_{8}$}
\label{J8}
\end{center}
\end{figure}

Let $B_{31}$ be the $(-2)$-curve $\tilde
L_{00,13,22}:(\l_{2}-1)x_{1}+x_{2}-\l_{2}=0$. Then $B_{32}$ and
$B_{31}$ form a fiber of type I$_{2}$ at $u=1$. Also $A_{32}$ and the
pullback of a certain $(2,2)$-curve form another fiber of type I$_{2}$
at $u=1/(\l_{1}\l_{2})$, while $A_{31}$ together with the pullback of
a certain $(2,2)$-curve form the third fiber of type I$_{2}$ at
$u=(\l_{1}-1)^{-1}(\l_{2}-1)^{-1}$. The change of variables
\begin{align*}
&X =u\bigl((\l_{1}-1)(\l_{2}-1)u-1\bigr)
\frac{(x_{2}-1)( \l_{2}x_{1}-x_{2})}{(\l_{2}-1) x_{2} (x_{1}-1)},
\\
&Y = -u^{3}\bigl((\l_{1}-1)(\l_{2}-1)u-1\bigr)
\frac{\l_{2} (x_{2}-1)( \l_{2}x_{1}-x_{2})}{t\,x_{2} (x_{1}-1)},
\end{align*}
converts \eqref{Kummer} to the Weierstrass equation
\begin{multline*}
Y^2 = X^3-u\bigl((2\l_{1}\l_{2}-\l_{1}-\l_{2}+2)u-2\bigr)X^2
\\
-u^2(u-1)(\l_{1}\l_{2}u-1)\bigl((\l_{1}-1)(\l_{2}-1)u-1\bigr)X.
\end{multline*}
Its discriminant is
\begin{multline*}
\Delta(u)=16u^{8}(u-1)^{2}(\l_{1}\l_{2}u-1)^{2}\bigl((\l_{1}-1)(\l_{2}-1)u-1\bigr)^{2}\\
\times \bigl(4\l_{1}\l_{2}(\l_{1}-1)(\l_{2}-1)u+(\l_{1}-\l_{2})^2\bigr).
\end{multline*}
[If $\l_{2}=-\l_{1}, 2-\l_{2}$, or $\l_{1}/(2\l_{1}-1)$, this elliptic
fibration has fiber of type III for general $\l_{1}$.]

Generically, it has only one section other than $0$-section:
\[
G_{2}  \quad\leftrightarrow\quad  T = (0,0).
\]

\subsection{$\J_{11}$}
Modifying the divisors appearing in the type $\J_{7}$ fibration we
constructed in \S4.2, we form two divisors
\begin{align*}
&\Psi_{11,0} = A_{31} + A_{21} + 2 G_{1}+ 2 A_{01} + 2 F_{0} + 2 A_{03} + 2 G_{3} + A_{33} + B_{33},
\\
&\Psi_{11,\infty} =  A_{30} + A_{20} + 2 G_{0}+ 2 A_{10} + 2 F_{1} + 2 A_{12} + 2 G_{2} + A_{32} + B_{32}.\end{align*}
They are of type I$_{4}^{*}$ and they do not intersect with each
other.
\begin{figure}[!ht]
\begin{center}
\[
\includegraphics[scale=0.72]{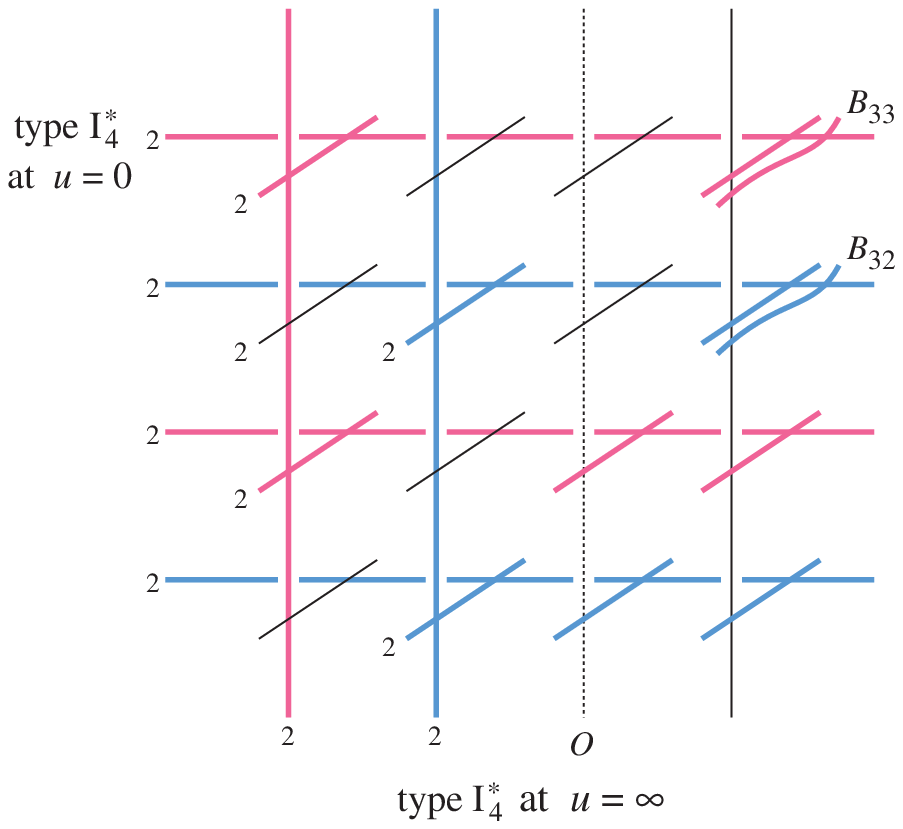}
\]
\caption{$\J_{11}$}
\label{J11}
\end{center}
\end{figure}

We look for a function whose divisor is $\Psi_{11,0} -
\Psi_{11,\infty}$, and we obtain the elliptic parameter of type
$\J_{11}$:
\begin{equation*}\label{ell_param:J_11}
u =
\frac{x_{2}(x_{2}-\l_{2})(x_{1}-x_{2})}{x_{1}(x_{2}-1)(\l_{2}x_{1}-x_{2})}.
\end{equation*}
The change of variables
\begin{align*}
X&= u \frac{(\l_{1}-1)(x_{2}-\l_{2})(x_{1}-x_{2})}{x_{1}(x_{1}-1)} \\
&=\frac{(\l_{1}-1)\,x_{2}(x_{2}-\l_{2})^{2}(x_{1}-x_{2})^{2}}
{x_{1}^{2}(x_{1}-1)(x_{2}-1)(\l_{2}x_{1}-x_{2})},
\\
Y&=
u^{2}\frac{(\l_{1}-1)(x_{2}-\l_{2})^{2}(x_{1}-x_{2})^{2}}{t\,x_{1}^{2}(x_{1}-1)^{2}} \\
&=\frac{(\l_{1}-1)\,x_{2}^{2}(x_{2}-\l_{2})^{4}(x_{1}-x_{2})^{4}}
{t\,x_{1}^{4}(x_{1}-1)^{2}(x_{2}-1)^{2}(\l_{2}x_{1}-x_{2})^{2}}
\end{align*}
converts \eqref{Kummer} to the Weierstrass equation
\begin{multline*}
Y^{2} = X^{3} +\bigl(\l_{1}u^2 - (2\l_{1}\l_{2}-\l_{1}-
\l_{2}+2)u+\l_{2}\bigr)\,u\,X^{2}
\\
+ (\l_{1}-1)(\l_{2}-1)\bigl((\l_{1}\l_{2}+1)u-2\l_{2}\bigr)\,u^{3}\,X +
\l_{2}(\l_{1}-1)^{2}(\l_{2}-1)^2\,u^{5}.
\end{multline*}
Its discriminant is of the form $u^{10}d(u)$, where $d(u)$ is a
polynomial of degree~$4$. The discriminant of $d(u)$ is too
complicated to write down here. However, a simple search reveals that
there are cases where four I$_{1}$ fibers degenerate even when $C_{1}$
and $C_{2}$ are not isogenous.
\begin{remark}
Suppose that the characteristic of the base field is~$0$.

(1) \ If $\l_{1}=-1$ and $\l_{2}=9\pm4\sqrt{5}$, then the fibration
has one I$_{2}$ fiber and one type II fiber. In this case
$j$-invariant of $C_{1}$ is $1728$ and that of $C_{2}$ is
$78608=2^{4}17^{3}$. They are not isogenous, and they can be defined
over~$\Q$.

(2) \ If $\l_{1}=-1$ and $\l_{2}=\pm\sqrt{-1}$, then the fibration has
two type II fibers. In this case $j$-invariant of $C_{1}$ is $1728$ and
that of $C_{2}$ is $128$. They are not isogenous, and they can be
defined over~$\Q$.
\end{remark}

\section{$(2,2)$-curves and $\J_{5}$, $\J_{9}$ and $\J_{10}$}

\subsection{$(2,2)$-curves}
Now the pullbacks of $(1,1)$-curves are not enough to construct all
the elliptic fibrations in Oguiso's list. A pullback of a
$(2,2)$-curve is a candidate for missing $(-2)$-curves. A nonsingular
$(2,2)$-curve in $\Proj^{1}\times\Proj^{1}$ is a curve of genus~$1$,
and thus, we first look for $(2,2)$-curves with a node. Then we try to
impose conditions such that their pullbacks are $(-2)$-curves. Here,
we do not try to make a systematic search as before.

Actually, we can construct an elliptic fibration of type $\J_{5}$ 
using only pullbacks of $(1,1)$-curves and the basic curves.
As a by-product, however, we obtain some new $(-2)$-curves 
which are pullbacks of $(2,2)$-curves. Such curves
have a node at $R_{11}$. They are given by an equation of the form
\[
a\, x_{1}^2z_{2}^2 + b\,x_{2}x_{1}z_{2}z_{1}+c\,x_{2}^2z_{1}^2+d\,x_{1}z_{2}^2z_{1}
+e\,x_{2}z_{2}z_{1}^2+f\,z_{2}^2z_{1}^2 = 0. 
\]
The fact that it has a node at $R_{11}$ corresponds to the fact that
the equation does not have the terms $x_{1}^{2}x_{2}^{2}$,
$x_{1}^{2}x_{2}z_{2}$, and $x_{1}z_{1}x_{2}^{2}$. In order to obtain
such a $(2,2)$-curve, we need to specify six points among $R_{ij}$
($1\le i,j \le 4$) such that no three among them are on the same
$F_{i}$ or $G_{j}$. We use such curves to construct an elliptic
fibration of type $\J_{9}$ and~$\J_{10}$.

\subsection{$\J_{5}$}

An elliptic fibration of type $\J_{5}$ has six I$_{2}$ fibers together
with one  I$_{6}^{*}$ fiber. In order to write down an elliptic
parameter for $\J_{5}$, we need to identify these six I$_{2}$ fibers.

Let $B_{33}$ and $B_{32}$ be the $(-2)$-curves introduced in \S\ref{notation}.
Consider two more $(-2)$-curves of this type:
\begin{align*}
&B_{12}:= \tilde L_{00,23,31} : \l_{2}(x_{1}-\l_{1})+(\l_{1}-1)x_{2}=0, \\
&B_{13}:= \tilde L_{00,22,31}: x_{1} - \l_{1} + (\l_{1}- 1) x_{2} = 0.
\end{align*}

\begin{figure}[!ht]
\begin{center}
\[
\includegraphics[scale=0.8]{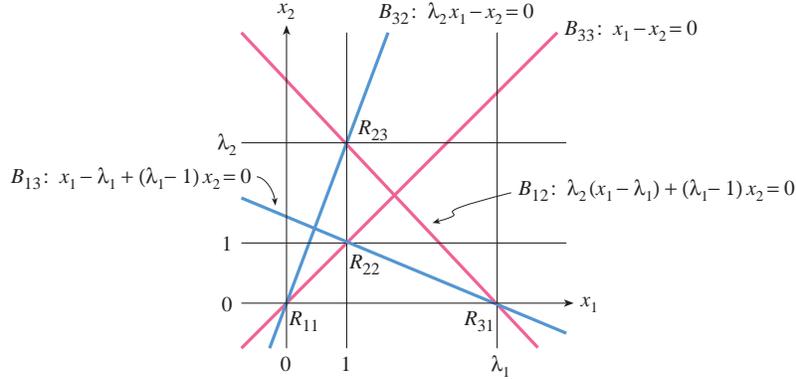}
\]
\caption{$(1,1)$-curves}
\label{J5_plane}
\end{center}
\end{figure}

Looking at Fig.~\ref{J5_plane}, we see that $B_{33}$ and $B_{12}$
intersect each other only at two points above the intersection of
lines $x_{1}-x_{2}=0$ and $\l_{2}(x_{1}-\l_{1})+(\l_{1}-1)x_{2}=0$.
Thus, the divisor $B_{33}+B_{12}$ is a singular fiber of type I$_{2}$.
Similarly, $B_{32}+B_{13}$ is another singular fiber of type I$_{2}$.
Furthermore, $B_{33}+B_{12}$ and $B_{32}+B_{13}$ do not intersect each
other since the image of these curves in $\A_{x_{1}}\times \A_{x_{2}}$
intersect only at $R_{ij}$ (see Fig.~\ref{J5_plane} below).

Computing the divisors $(x_{1}-x_{2})$,
$\bigl(\l_{2}(x_{1}-\l_{1})+(\l_{1}-1)x_{2}\bigr)$,
$(\l_{2}x_{1}-x_{2})$ and $\bigl(x_{1}-\l_{1}+(\l_{1}-1)x_{2}\bigr)$,
we see that
\[
u = \frac{(x_{1}-x_{2})\bigl(\l_{2}(x_{1}-\l_{1})+(\l_{1}-1)x_{2}\bigr)}
{(\l_{2}x_{1}-x_{2})\bigl(x_{1}-\l_{1}+(\l_{1}-1)x_{2}\bigr)}
\]
is an elliptic parameter of type $\J_{5}$.  We have
\[
u-1 = \frac{-\l_{1} (\l_{2}-1)\,x_{2} (x_{1}-1)}
{(\l_{2}x_{1}-x_{2})\bigl(x_{1}-\l_{1}+(\l_{1}-1)x_{2}\bigr)},
\]
which shows that the fiber at $u=1$ is a singular fiber of type
I$_{6}^{*}$.
\begin{figure}[!ht]
\begin{center}
\[
\includegraphics[scale=0.72]{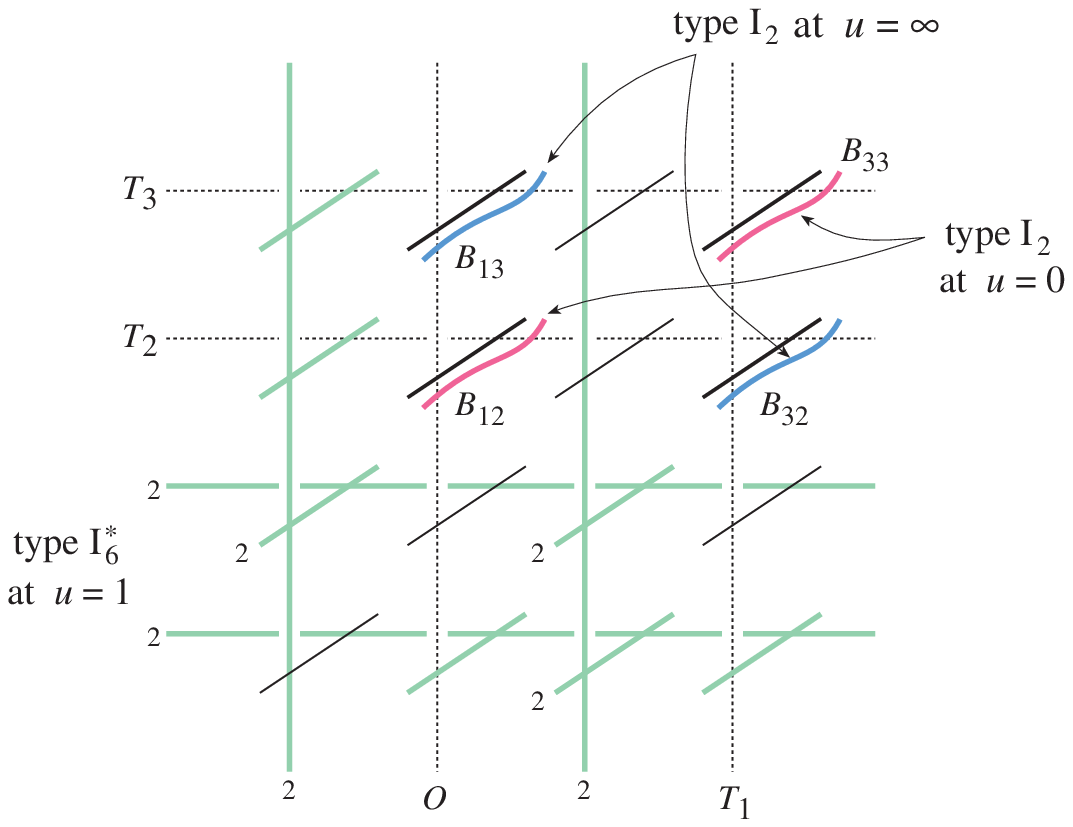}
\]
\caption{$\J_{5}$}
\label{J5}
\end{center}
\end{figure}
Each of the divisors $A_{12}$, $A_{13}$, $A_{32}$ and $A_{33}$ is a
component of a singular fiber of type I$_{2}$. The other $(-2)$-curves
are pullbacks of $(2,2)$-curves. For example, the singular fiber at
$u=\l_{1}\l_{2}-\l_{1}+1$ consists of $A_{12}$ and the pullback of the
$(2,2)$-curve given by
\[
\l_{2}x_{1} ( x_{1}-1) + (\l_{1}-1)(x_{2}-1)(\l_{2}x_{1}-x_{2}) =0.
\]

In order to obtain a Weierstrass equation using the curve $G_{0}$ as
the $0$-section, we first put
\begin{align*}
X_{0}
&=\frac{(x_{1}-x_{2})(x_{1}-\l_{1})}{x_{1}\bigr((x_{1}-\l_{1})+(\l_{1}-1)x_{2}\bigr)},
\\
Y_{0} &= \frac{\l_{1}(\l_{1}-1)
\, t\, x_{2}(x_{1}-1)(x_{1}-\l_{1})(x_{2}-x_{1})}
{x_{1}(\l_{2}x_{1}-x_{2})\bigr((x_{1}-\l_{1})+(\l_{1}-1)x_{2}\bigr)^2},
\end{align*}
and then put
\begin{align*}
X &= \l_{1} (\l_{2}-1)(u-1)(u-\l_{1}\l_{2}+\l_{1}-1)
\bigl((\l_{1}\l_{2}-\l_{1}-\l_{2})u+\l_{2}\bigr) X_{0},
\\
Y &= \l_{1}^2(\l_{2}-1)^2(u-1)^2(u-\l_{1}\l_{2}+\l_{1}-1)
\bigl((\l_{1}\l_{2}-\l_{1}-\l_{2})u+\l_{2}\bigr)Y_{0}.
\end{align*}
Then $(X,Y)$ satisfy the Weierstrass equation
\[
Y^{2}= X(X-\alpha)(X-\beta),
\]
where
\begin{align*}
\alpha &= -\l_{1} (\l_{2}-1) (u-1)\bigl((\l_{1}\l_{2}-1)u-\l_{1}+1\bigr)
\bigl((\l_{1}\l_{2}-\l_{1}-\l_{2})u+\l_{2}\bigr),
\\
\beta &=\l_{1} (\l_{2}-1) u(u-1)(u-\l_{1}\l_{2}+\l_{1}-1)
\bigl((\l_{1}\l_{2}-\l_{2})u-\l_{1}+\l_{2}\bigr).
\end{align*}
The discriminant of this fibration is given by
\begin{multline*}
\Delta(u) =16\l_{1}^{6}\l_{2}^{2}(\l_{1}-1)^{2}u^{2}(u-1)^{12}
\\
\times
(u-\l_{1}\l_{2}+\l_{1}-1)^{2}\bigl((\l_{1}\l_{2}-1)u-\l_{1}+1\bigr)^{2}
\\
\times
\bigl((\l_{1}\l_{2}-\l_{2})u-\l_{1}+\l_{2}\bigr)^{2}
\bigl((\l_{1}\l_{2}-\l_{1}-\l_{2})u+\l_{2}\bigr)^{2}.
\end{multline*}

The Mordell-Weil group of this elliptic surface has the following
three sections:
\[\renewcommand{\arraystretch}{1.41421356}
\begin{array}{ccl}
F_{3} & \leftrightarrow & T_{1} = (0,0),
\\
G_{2} & \leftrightarrow & T_{2} = (\alpha,0),
\\
G_{3} & \leftrightarrow & T_{3} = (\beta, 0).
\end{array}
\]

\subsection{$\J_{9}$}

In order to construct an elliptic fibration of type $\J_{9}$, we need
to find a divisor of type I$_{0}^{*}$ different from the ones
appearing in $\J_{4}$ or $\J_{7}$. To do so we look for a $(-2)$-curve
$P_{33}$ such that $2G_{3} + A_{03} + A_{33} + B_{33} +P_{33}$ is of
type I$_{0}^{*}$. We can show that $P_{33}$ cannot be a pullback of a
$(1,1)$-curve; if that were the case, $B_{33}$ and $P_{33}$ would have
to intersect each other. Thus, we look for a $(2,2)$-curve whose
double cover serves as $P_{33}$.
\begin{figure}[!ht]
\begin{center}
\[
\includegraphics[scale=0.8]{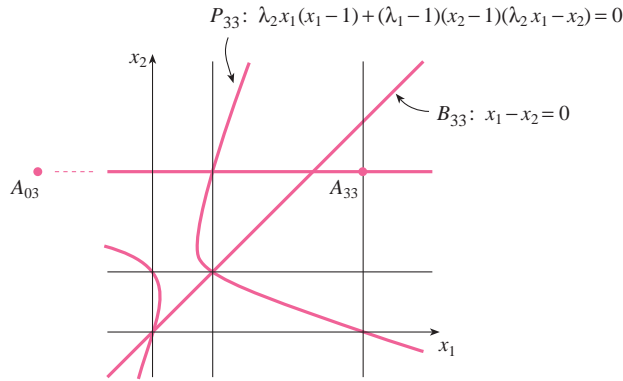}
\]
\caption{fiber at $u=0$}
\label{J9_at_0}
\end{center}
\end{figure}

It turns out that the pullback of the $(2,2)$-curve
\begin{equation}\label{(2,2)-curve} \l_{2}x_{1} (
x_{1}-1)+(\l_{1}-1)(x_{2}-1)(\l_{2}x_{1}-x_{2})=0 \end{equation} can
be used as $P_{33}$. This curve is a component of a I$_{2}$ fiber of
the elliptic fibration of type $\J_{5}$ which we constructed in the
previous subsection. The $(2,2)$-curve \eqref{(2,2)-curve} has a node
at $R_{00}$, and passes through $R_{11}, R_{12},R_{22},R_{23}$, and
$R_{31}$. Fig.~\ref{J9_at_0} shows the projection of the $(-2)$-curves
contained in the divisor $\Psi_{9,0} = 2G_{3} + A_{03} + A_{33} +
B_{33} +P_{33}$. (The projection of $A_{03}$ is $R_{03}$, which is a
point at infinity.)

\begin{figure}[!ht]
\begin{center}
\[
\includegraphics[scale=0.8]{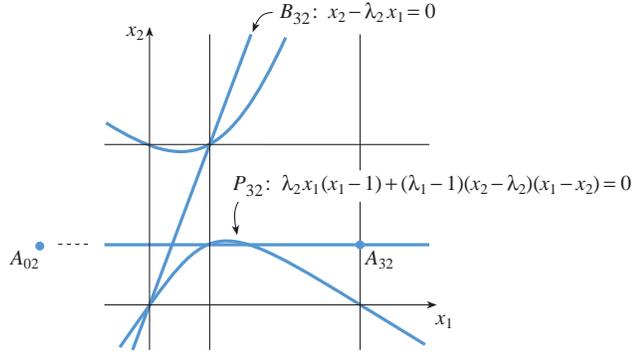}
\]
\caption{fiber at $u=\infty$}
\label{J9_at_infty}
\end{center}
\end{figure}

\begin{figure}[!ht]
\begin{center}
\[
\includegraphics[scale=0.72]{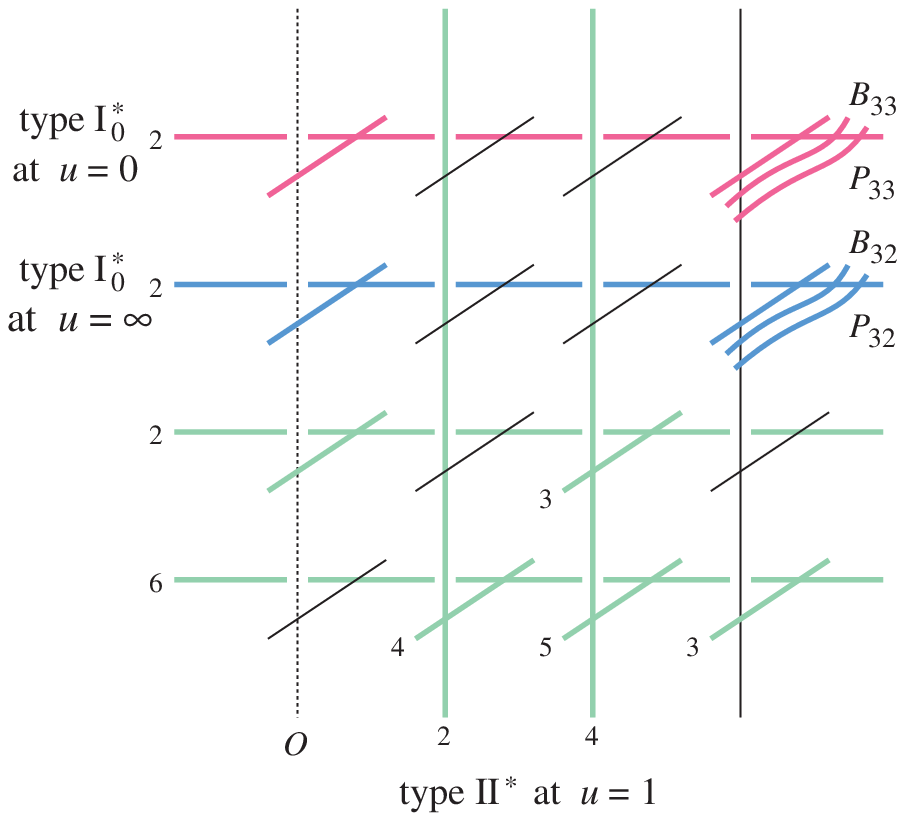}
\]
\caption{$\J_{9}$}
\label{J9}
\end{center}
\end{figure}

Similarly, let $P_{32}$ be the pullback of the $(2,2)$-curve
\[
\l_{2}x_{1} ( x_{1}-1 )+(\l_{1}-1)( x_{2}-\l_{2} ) (x_{1}-x_{2})=0.
\]
Then, the divisor $\Psi_{9,\infty} = 2G_{2} + A_{02} + A_{32} + B_{32}
+P_{32}$ is again of type I$_{0}^{*}$, which does not intersect with
$\Psi_{9,0}$. Fig.~\ref{J9_at_infty} shows the curves contained in the
divisor $\Psi_{9,\infty}$.
Looking for the function having the divisor
$\Psi_{9,0}-\Psi_{9,\infty}$, we find the elliptic parameter $u$ given
by
\begin{equation}\label{ell_param:J_9}
u = \frac{(x_{2}-\l_{2})(x_{1}-x_{2}) 
\bigl(\l_{2}x_{1} ( x_{1}-1)+(\l_{1}-1)(x_{2}-1)(\l_{2}x_{1}-x_{2})\bigr)}
{(x_{2}-1)(\l_{2}x_{1} -x_{2}) 
\bigl(\l_{2}x_{1} ( x_{1}-1 )+(\l_{1}-1)( x_{2}-\l_{2} )
(x_{1}-x_{2})\bigr)}.
\end{equation}
We have
\begin{multline*}
u -1 = \\
\frac{-\l_{2}(\l_{2}-1)x_{1}x_{2}(x_{1}-1)^{2}}
{(x_{2}-1)(\l_{2}x_{1} -x_{2}) 
\bigl(\l_{2}x_{1} ( x_{1}-1 )+(\l_{1}-1)( x_{2}-\l_{2} )(x_{1}-x_{2})\bigr)}.
\end{multline*}
The zero divisor of this function $u-1$ is given by the following
divisor consisting only of the basic curves:
\[
A_{01} + 2 G_{1} + 3 A_{21} + 4 F_{2} + 5 A_{20} + 6 G_{0} + 3 A_{30}
+ 4 A_{10} + 2 F_{1}.
\]
This is a singular fiber of type II$^{*}$ (see Fig.~\ref{J9}). Thus,
the elliptic parameter given by~\eqref{ell_param:J_9} is of type
$\J_{9}$.

Our next task is to write down a Weierstrass equation. If we regard
\eqref{ell_param:J_9} as the defining equation of a curve in
$\Proj^{1}\times \Proj^{1}$ defined over $k(u)$, then we can show that
this curve is a curve of genus~$0$, and thus it can be
parametrized. In fact, we can parametrize $x_{1}$ and $x_{2}$
satisfying \eqref{ell_param:J_9} using the parameter
\[
\xi=-\frac{(x_{2}-1)(\l_{2}x_{1} - x_{2})}{\l_{2} x_{1}x_{2}}.
\]
Actual parametrizations of $x_{1}$ and $x_{2}$ are complicated and we
omit here.  Substituting $x_{1}$ and $x_{2}$ in the equation
\eqref{Kummer} by these parametrizations, we obtain an equation of a
curve of genus~$1$ with variables in $(\xi, t)$ defined over $k(u)$.
This equation turns out to be a quadratic equation in $t$, and it is
easily converted to a Weierstrass equation.
Combining all these, we obtain the change of variables
\[
X_{0} =-\frac{\l_{2}(\l_{2}-1)
x_{1}(x_{1}-1)}{(x_{2}-1)(\l_{2}x_{1}-x_{2})}, \quad Y_{0}
=\frac{\l_{2}(\l_{2}-1)}{t}
\]
that converts \eqref{Kummer} to a twisted form of Weierstrass equation
\begin{multline*}
u(u-1)Y_{0}^{2} = X_{0}^{3} + (\l_1 \l_2 -2 \l_1-2\l_2+1) (u-1) X_{0}^2 \\
- (\l_1 + \l_2 -1)(\l_1 \l_2 - \l_1-\l_{2}) (u-1)^2 X_{0} \\
- \l_1 \l_2 (\l_1 - 1)(\l_2 - 1) (u-1)^2.
\end{multline*}
By letting
\[
X =u(u-1)X_{0}, \quad Y= u^{2}(u-1)^{2}Y_{0},
\]
we obtain the following Weierstrass equation:
\begin{multline*}
Y^{2} = X^{3} + (\l_1 \l_2 -2 \l_1 -2\l_2+1)\,u\,(u-1)^2 X^2 \\
- (\l_1 + \l_2 -1)(\l_1 \l_2 - \l_1-\l_{2})\,u^2 (u-1)^4 X \\
- \l_1 \l_2(\l_1 - 1) (\l_2 - 1)\,u^3 (u-1)^5.
\end{multline*}
Its discriminant is of the form $u^{6}(u-1)^{10}d(u)$, where $d(u)$ is
a polynomial of degree~$2$. The discriminant of $d(u)$ is given by
\[
16 \l_{1}^{2}\l_{2}^{2}(\l_{1}-1)^2(\l_{2}-1)^{2}
(\l_{1}^2-\l_{1}+1)^3(\l_{2}^2-\l_{2}+1)^3.
\]
If either $\l_{1}$ or $\l_{2}$ is a sixth root of unity, then two
I$_{1}$ fibers of the fibration degenerate to form a type II fiber.

\subsection{$\J_{10}$}
\begin{figure}[!ht]
\begin{center}
\[
\includegraphics[scale=0.72]{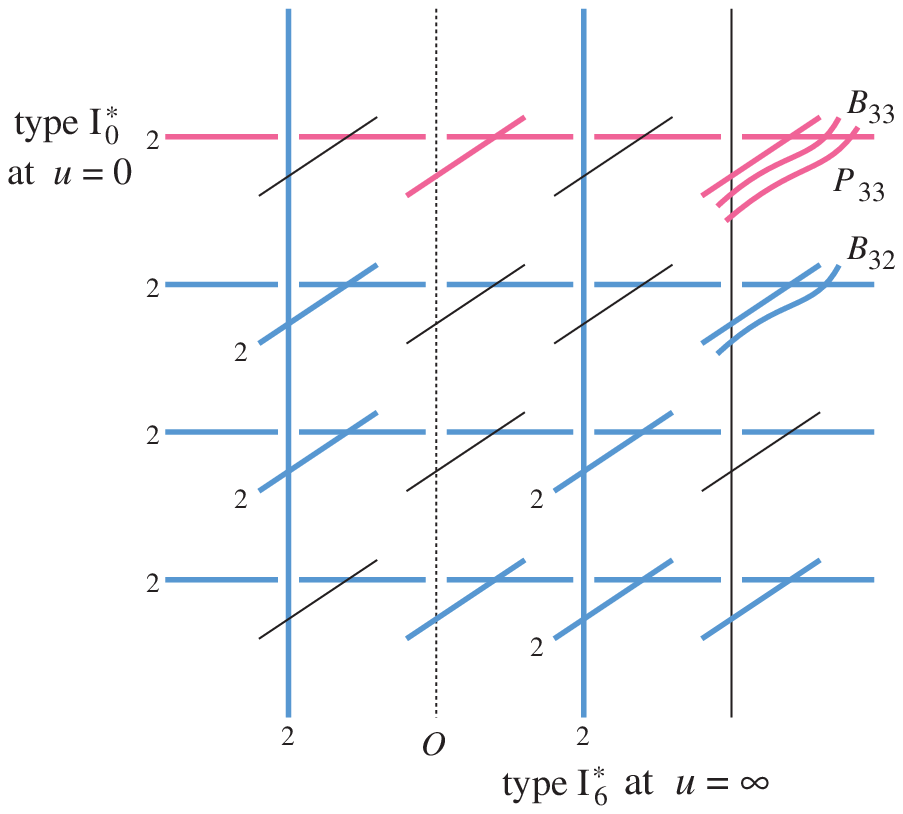}
\]
\caption{$\J_{10}$}
\label{J10}
\end{center}
\end{figure}

In order to construct an elliptic fibration of type $\J_{10}$, we must
find yet another divisor of type I$_{0}^{*}$. The divisor $\Psi_{9,0}
=2 G_{3} + A_{03} + A_{33} + B_{33} + P_{33}$ is a divisor of type
I$_{0}^{*}$ appearing in the elliptic fibration constructed in the
previous subsection. Since neither $B_{33}$ nor $P_{33}$ intersects
with $A_{13}$, we see that
\[
\Psi_{10,0} = 2 G_{3} + A_{13} + A_{33} + B_{33} + P_{33}
\]
is also a divisor of type I$_{0}^{*}$. We then find a divisor of
type I$_{6}^{*}$ that does not intersect with $\Psi_{10,0}$:
\begin{multline*}
\Psi_{10,\infty} = B_{32} + A_{32} + 2 G_{2} + 2 A_{02} + 2 F_{0} + 2 A_{01} \\
+ 2 G_{1} + 2 A_{21} + 2 F_{2} + 2 A_{20} + 2 G_{0} + A_{10} + A_{30}
\end{multline*}
(see Fig.\ref{J10}). Looking for the function having the divisor
$\Psi_{10,0}-\Psi_{10,\infty}$, we find the elliptic parameter of type
$\J_{10}$ given by
\begin{equation}\label{ell_param:J_10}
u = \frac{( x_{2}-\l_{2} )(x_{1}-x_{2})
\bigl((\l_{1}-1)( x_{2}-1 )(\l_{2}x_{1} -x_{2}) +\l_{2}x_{1} ( x_{1}-1)\bigr)}
{x_{2}(x_{2}-1)(x_{1}-1)(\l_{2}x_{1} -x_{2})}.
\end{equation}
The curve in $\Proj^{1}\times \Proj^{1}$ over $k(u)$ defined by
\eqref{ell_param:J_10} is a curve of genus~$0$.  As in the
case of $\J_{9}$, the parameter
\[
\xi=\frac{(x_{1}-x_{2})(x_{2}-\l_{2})}{(x_{1}-1)x_{2}},
\]
can be used to parametrize this curve. We can proceed in a similar
manner to the case of $\J_{9}$ and we obtain the change of
variables
\begin{align*}
&X_{0}=\frac{\l_{1}(\l_{1}-1)(x_{2}-1)(\l_{2}x_{1}-x_{2})}{x_{1}(x_{1}-1)},
\\
&Y_{0}=
\frac{\l_{1}(\l_{1}-1)(x_{2}-1)^2(\l_{2}x_{1}-x_{2})^2}{t x_{1}^2(x_{1}-1)^2},
\end{align*}
which converts \eqref{Kummer} to
\begin{multline*}
uY_{0}^{2} = X_{0}^{3} - (u +\l_{1}+\l_{2}-1)(u-\l_{1}\l_{2}+\l_{1}+\l_{2}) 
X_{0}^{2}
\\
+\l_{1}\l_{2}(\l_{1}-1)(\l_{2}-1)(2u-\l_{1}\l_{2}+2\l_{1}+2\l_{2}-1)X_{0}
\\
+ \l_{1}^{2}\l_{2}^{2}(\l_{1}-1)^2 (\l_{1}-1)^2.
\end{multline*}
Putting
\[
X =u X_{0}, \quad Y = u^{2} Y_{0},
\]
we obtain the Weierstrass equation
\begin{multline*}
Y^{2} = X^{3} + u(u-\l_{1}-\l_{2}+1)(u +\l_{1}\l_{2}-\l_{1}-\l_{2}) X^{2}
\\
+\l_{1}\l_{2}(\l_{1}-1)(\l_{2}-1)u^2(2u+\l_{1}\l_{2}-2\l_{1}-2\l_{2}+1)X
\\
+ \l_{1}^{2}\l_{2}^{2}(\l_{1}-1)^2 (\l_{1}-1)^2u^3.
\end{multline*}
Its discriminant is of the form $u^{6}d(u)$, where $d(u)$ is a
polynomial of degree~$2$. We can show that $d(u)$ can have a multiple
root without $C_{1}$ and $C_{2}$ being isogenous.

\section{Full list of the defining equations in a special case}  

In this section, we take as $C_{1}$ and $C_{2}$  the most familiar elliptic curves 
\begin{equation}
C_1: y_1^2=x_1^3 - x_1,  \quad C_2: y_2^2=x_2^3 - 1,
\end{equation}
and write down the full list of the defining equations of mutually nonisomorphic elliptic fibrations on the Kummer surface $S=\Km (C_1 \times C_2)$ in characteristic~$0$.
Although they  are very special among elliptic curves (e.g. automorphisms or complex multiplications), corrseponding Kummer surface $S$ serves as a more or less 
``typical" case, since $C_{1}, C_{2}$ are not isogenous to each other. 

In this case,   the number $N(n)$ of  nonisomorphic elliptic fibrations on $S$ of type $\J_{n}$ has been determined by Oguiso as follows:
$$ N(n)= 1 \quad \mbox{for} \  n = 2, 3, 5, 8, 9, 10, $$
and
$$ N(n)=2 \quad \mbox{for} \    n=1, 4, 6, 7, 11.$$
(See Oguiso \cite[p.~652]{O}. We note that this number $N$ is {\em not} typical among all non-isogenous curves, as shown there.)

Now observe that the values of Legendre parameter $\lambda_i$ for the present $C_i$ are as follows:
$$ \lambda_1= -1, -2 \ \mbox{or}\ 1/2, \quad \lambda_2=-\omega \ \mbox{or}\         -\omega^2,
$$
where $\omega$ is a cubic root of unity.
In the following, we write down the $N=N(n)$ defining equations for each type
 $\J_{n}$. When $N=1$,  we give essentially the same equation as the one constructed in the previous sections, except that we make some coordinate change when it makes the equation look simpler.
   When $N>1$, we make the  same construction as  before using a suitable equivalent value of $\lambda_i$. We briefly indicate how to verify that the resulting defining equations are not isomorphic to each other.

\subsection{$\J_{1}$}
\begin{equation} \label{special:J_1(i)}
 y^2=x \bigl(x^2+(u^4+1) x + 4 u^4\bigr), 
\end{equation}
$$
 J=\frac{1}{108} \frac{(u^{8}-10 u^4+1)^3}{u^8 (u^8-14u^4+1)}.
$$

\begin{equation}\label{special:J_1(ii)}
  y^2=x\bigl(x^2+(u^{4}+6 (2\omega+1) u^2+1) x -32 u^4\bigr), 
\end{equation}
$$
 J=\frac{1}{6912}
 \frac{\bigl(u^8+12 (2\omega+1) u^6-10 u^4+12 (2\omega+1) u^2+ 1\bigr)^3}
{u^8 \bigl(u^8+12 (2\omega+1) u^6 + 22 u^4+12 (2\omega+1) u^2+ 1\bigr)}.
$$
Both the equations (\ref{special:J_1(i)}) and (\ref{special:J_1(ii)}) 
 have two I$_8$ fibers at $u=0$ and $\infty$ and eight I$_1$ fibers.
Suppose they define isomorphic elliptic curves over $k(u)$. Then there must be a linear transformation of $u$ fixing $0$ and $\infty$ which sends one $J$ into the other, $J$ denoting the classical absolute invariant of the generic fibre (normalized so that $J=1$ for $y^2=x^3-x$).  But this is impossible, as the positions of the eight  I$_1$ fibers are determined by the simple poles of $J$ and they cannot be transformed by such a linear transformation. This proves that the two elliptic fibrations are not isomorphic to each other. 

\subsection{$\J_{2}$}
\begin{equation} \label{special:J_2}
 y^2=x\bigl(x^2-(3 u^4+6 u^{2}-1) x + 32 u^6\bigr), 
\end{equation}
$$
 J=\frac{1}{6912} \frac{(9 u^8 -60 u^6 + 30 u^4 -12 u^2  + 1)^3}
 {u^{12} (u^4-10 u^2+1)(9 u^4-2u^2+1)}.
$$

\subsection{$\J_{3}$}
\begin{equation} \label{special:J_3}
 y^2=x^3 + u^4 (u^4+1),    \quad J=0.
\end{equation}

\subsection{$\J_{4}$}
\begin{equation} \label{special:J_4(i)}
 y^2=x^3 - (u^3-1)^2 x, \ (u=x_1), \quad J=1.
\end{equation}
\begin{equation} \label{special:J_4(ii)}
 y^2=x^3 - (v^3-v)^3, \ (v=x_2), \quad J=0.
\end{equation}
These are the two obvious elliptic fibrations on $S$ induced by the 
projections 
$C_1 \times C_2 \to C_1$ or $C_2$.

\subsection{$\J_{5}$}
\begin{equation} \label{special:J_5}
 y^2=x(x-4)\bigl(x+2u(u^2+3 u+3)\bigr).
\end{equation}
$$
J=\frac{1}{27}\frac{(u^6+6 u^5 + 15 u^4 + 20 u^3 +15 u^2 +6 u  +4)^3}
{u^2(u+2)^2(u^2+u+1)^2 (u^2+3 u+3)^2}.
$$

\subsection{$\J_{6}$}
\begin{equation} \label{special:J_6(i)}
 y^2=x(x+2 u^2) \bigl(x -  u (u^2-u+1)\bigr), 
\end{equation}
$$
 J=\frac{1}{27} \frac{(u^4+5 u^2+1)^3}{u^2 (u^4+u^2+1)^2}.
$$
\begin{equation} \label{special:J_6(ii)}
 y^2=x (x-\omega u^2)\bigl(x + u(2u-1) (u+\omega^2)\bigr), 
\end{equation}
$$
 J=\frac{4}{27} 
 \frac{\omega\bigl(2u^{2}-(\omega +2)u-\omega^2\bigr)^3
 \bigl(2u^{2}-2(\omega +2) u-\omega^2\bigr)^3}%
 {u^2(u-1)^2 (2u-1)^2(u+\omega^{2})^2 (2 u+\omega^{2})^2}.
$$
The Legendre parameters we employed for the first equation (\ref{special:J_6(i)}) are $\lambda_1=-\omega, \lambda_2 =-1$, while those for the second one (\ref{special:J_6(ii)}) are $\lambda_1=-\omega, \lambda_2 =2$. That the two equations define nonisomorphic elliptic fibrations  
 can be checked in the same way as the case for $\J_{1}$ above.

\subsection{$\J_{7}$}
\begin{equation} \label{special:J_7(i)}
 y^2=x\bigl(x^2-u(u+1)(u+2) x + u^2(u+2)^2\bigr), 
\end{equation}
$$
 J=\frac{4}{27} \frac{(u^2+2u-2)^3}{(u-1)(u+3)}.
$$

\begin{equation} \label{special:J_7(ii)}
  y^2=x\bigl(x^2+ \omega u(u-1)(u-3\omega-2) x +2\omega^2 u^2(u-1)^2\bigr).
\end{equation}
$$
 J=\frac{1}{27} \frac{(u^2-2(3\omega+2)u+3\omega-11)^3}{(u^2-2(3\omega+2)u+3\omega-13)}.
$$
In this case, we can check  that
 there is a linear transformation 
of $u$ sending the first $J$ into the second one. However
it does not preserve the position of singular fibres which 
can be seen from the discriminants (but not from the 
absolute invariants). Hence 
(\ref{special:J_7(i)}) and (\ref{special:J_7(ii)}) are
not isomorphic.

\subsection{$\J_{8}$}
\begin{equation} \label{special:J_8}
 y^2=x\bigl(x^2+  u (3 u+2) x + u^2 (1+ 3u + 3u^2 + 2 u^3)\bigr), 
\end{equation}
$$
 J=\frac{4}{27} \frac{(6 u^3 - 3u -1)^3}{u^2 (2u+1)^2 (u^2+u+1)^2 (8u+3)}.
$$

\subsection{$\J_{9}$}
\begin{equation} \label{special:J_9}
 y^2=x^3+  u (u^2 -4)^3,  \quad J=0.
\end{equation}

\subsection{$\J_{10}$}
\begin{equation} \label{special:J_10}
 y^2=x^3+  u^2(u+3)x^2 +u^2(-2 u^2-2u+3) x +u^4(u-1).
\end{equation}
We omit $J$.

\subsection{$\J_{11}$}
\begin{multline} \label{special:J_11(i)}
 y^2=x^3-27u^2(u^4+6u^3+5u^2-6u+1) x \\
-54u^3(u^2+1)(u^4+9u^3+20u^2-9u+1).
\end{multline}
\begin{multline} \label{special:J_11(ii)}
 y^2=x^3- 27 u^2(4u^4-12u^3+10(\omega+1)u^2-6\omega u+\omega) x \\
+27 u^3\bigl(16  u^6-72 u^5+6(19+10 \omega)u^4 \\
-63(1+2\omega)u^3+3(-9+10\omega)u^2+18u-2\bigr). 
\end{multline}
We omit $J$, but it can be checked that the two elliptic fibrations are not isomorphic to each other by a similar argument as before.

Thus we have listed the defining equations of elliptic fibrations (with a section) on the Kummer surface $S=\Km (C_1 \times C_2)$ with $C_i$ given by (6.1) over an algebraically closed field $k$ of characteristic 0. Needless to say that  the function field $k(x,y,u)$ defined by each of the equations (6.2) through (6.17) is isomorphic to one and the same function field $k(S)$, which is the
extension $k(x_1,x_2,t)$ with
 $t=y_1/y_2$ determined by (6.1).

\section{Closing remark}
In closing this paper, it should be remarked that the problems posed
in the Introduction (\S 1.1) should be interesting and worth
considering for more general $K3$ surfaces.

Even in the case of Kummer surfaces, we could ask such questions as
follows:
\begin{problem} 
Study Problems 1 and 2 for the Kummer surface $X=\Km (A)$, when $A$ is
the Jacobian variety of a genus two curve.
\end{problem}
For this, the so-called $16_6$-configuration of thirty-two $(-2)$-curves
on $X$ should play an important role in place of the twenty-four basic
curves used in this paper. A special case has been treated in Shioda
\cite{ClKum}.

According to Weil \cite{Weil}, a principally polarized abelian surface
is either the Jacobian variety of a genus two curve or a product of
two elliptic curves. Beyond the case of principally polarized abelian
surfaces, we ask:
\begin{problem}  
Find at least one elliptic parameter for the Kummer surface $X=\Km
(A)$ when $A$ is a generic member in a family of polarized abelian
surfaces.
\end{problem}
The coefficients in the defining equation (especially the
discriminant) for such should be related to some modular forms or
theta-functions of interest.\\

{\bf Acknowledgements} \  
We would like to thank Professor Fumio Sakai of Saitama University who arranged a Seminar by the second author about ``Elliptic parameters'' in September 2005, which resulted in this fruitful joint work with the first author who attended the Seminar. We thank Matthias Sch\"utt for careful reading of the manuscript.
We are most grateful to the referee for making very sharp and constructive comments on the manuscript, who has requested the three things: (i) to make the paper free from mistakes (since many standard calculations are nicely omitted), (ii) to add a few words about the case of quasi-elliptic fibrations, and (iii) to include ``a full account in one non-trivial case'' after saying that Remark~1.7  is very impressive.

In answering these requests, (i) we have tried our best to minimize the possible errors (but still some errors might have crept in during correction). (ii) In $\operatorname{char}k=3$, it is possibe that some of the defining equations becomes a  quasi-elliptic fibration. We should leave this question to the interested reader,
but let us say this:
 if for example we let $\lambda_1=\lambda_2=-1$ in $\operatorname{char}k=3$, then $\J_3$-fibration  gives the same equation as the  $\J_3$-fibration in \S6 (for $\operatorname{char}k \neq 3$), which is indeed quasi-elliptic in $\operatorname{char}k=3$.
(iii)  We have  added a new section \S6 in the revised version to respond to this request. 

The second-named author ackowledges the  support from 
Grant-in-Aid for Scientific Research No.17540044.


\begin{thebibliography}{15}

\bibitem{Cassels}
J.~W.~S. Cassels,  Lectures on elliptic curves, London Mathematical
  Society Student Texts, \textbf{24}, Cambridge University Press, Cambridge, 1991.

\bibitem{Connell}
I.~Connell, Addendum to a paper of {K}. {H}arada and {M}.-{L}.\ {L}ang:
  ``{S}ome elliptic curves arising from the {L}eech lattice'' [{J}.\ {A}lgebra
 \textbf{125} (1989), no.\ 2, 298--310], J. Algebra \textbf{145} (1992),
  463--467.

\bibitem{Inose}
H.~Inose,  Defining equations of singular $K3$ surfaces and a notion of isogeny,
Proc. Intl. Symp. on Algebraic Geometry, Kyoto 1977, 
Kinokuniya, 1978, 495--502.

\bibitem{Kodaira}
K.~Kodaira, On compact analytic surfaces. {II}, {III}, Ann. of Math. (2)
  \textbf{77} (1963), 563--626; ibid. \textbf{78} (1963), 1--40.

\bibitem{Kuwata:can-height}
M.~Kuwata, The canonical height and elliptic surfaces, J. Number
  Theory \textbf{36} (1990), 201--211.

\bibitem{Kuwata:MWrank}
\bysame, Elliptic {$K3$} surfaces with given {M}ordell-{W}eil rank,
  Comment. Math. Univ. St. Paul. \textbf{49} (2000), 91--100.

\bibitem{Morrison}
D.~R. Morrison,  On {$K3$} surfaces with large {P}icard number, Invent.
  Math. \textbf{75} (1984), 105--121.

\bibitem{O}
K.~Oguiso,  On {J}acobian fibrations on the {K}ummer surfaces of the
  product of nonisogenous elliptic curves, J. Math. Soc. Japan \textbf{41}
  (1989), 651--680.

\bibitem{P-S-S}
I.~I. Pjatecki{\u\i}-{\v{S}}apiro  and I.~R. {\v{S}}afarevi{\v{c}},
   Torelli's theorem for algebraic surfaces of type ${\rm {K}}3$, Izv.
  Akad. Nauk SSSR Ser. Mat. \textbf{35} (1971), 530--572, English tansl.: Math.
  USSR, Izv. 5 (1971) 547--588.

\bibitem{IS}
T.~Shioda and H.~Inose, On singular {$K$3} surfaces, Complex Analysis
  and Algebraic Geometry, Iwanami Shoten, Tokyo, 1977, 119--136.

\bibitem{Shioda:MWL}
T.~Shioda, On the {M}ordell-{W}eil lattices, Comment. Math. Univ.
  St. Paul. \textbf{39} (1990), no.~2, 211--240.

\bibitem{K3SP}
\bysame,  A note on {$K3$} surfaces and sphere packings, Proc. Japan
  Acad. Ser. A Math. Sci. \textbf{76} (2000), 68--72.

\bibitem{ClKum}
\bysame,  Classical {K}ummer surfaces and {M}ordell-{W}eil lattices,
  Algebraic geometry, Contemp. Math., vol. 422, Amer. Math. Soc., Providence,
  RI, 2007, pp.~213--221.

\bibitem{Shioda:CorrMWL}
\bysame,  Correspondence of elliptic curves and {M}ordell-{W}eil lattices
  of certain elliptic {$K3$}'s, Algebraic Cycles and Motives, vol. 2, London
  Math. Soc. Lecture Note Ser., vol. 344, Cambridge Univ. Press, Cambridge,
  2007, pp.~319--339.

\bibitem{Weil}
A.~Weil,  Zum {B}eweis des {T}orellischen {S}atzes, Nachr. Akad.
  Wiss. G\"ottingen. Math.-Phys. Kl. IIa.  (1957), 33--53.

\end{thebibliography}
\providecommand{\bysame}{\leavevmode\hbox to3em{\hrulefill}\thinspace}

\end{document}